\documentclass[aap,MSNbibl,seceqn,citesort,dvips]{arximspdf}

\doi{10.1214/11-AAP771}
\volume{22}
\issue{3}
\pubyear{2012}
\firstpage{971}
\lastpage{1007}

\makeatletter

\newtheorem{Theorem}{Theorem}[section]
\newtheorem{Proposition}{Proposition}[section]
\newtheorem{Lemma}{Lemma}[section]

\newproclaim{Remark}{Remark}[section]

\newcommand{\ess}{\operatorname{ess}}
\newcommand{\esssup}{\ess\sup}

\newcommand{\doY}{\delta Y}
\newcommand{\dotY}{\delta\tY}
\newcommand{\doZ}{\delta Z}
\newcommand{\dof}{\delta f}

\newcommand{\eps}{\varepsilon}
\newcommand{\E}{\mathbb{E}}
\newcommand{\F}{\mathbb{F}}
\newcommand{\N}{\mathbb{N}}
\newcommand{\R}{\mathbb{R}}
\newcommand{\PP}{\mathbb{P}}

\newcommand{\Ac}{\mathcal{A}}
\newcommand{\Fc}{\mathcal{F}}
\newcommand{\Hc}{\mathcal{H}}
\newcommand{\Ic}{\mathcal{I}}
\newcommand{\Lc}{\mathcal{L}}
\newcommand{\Mc}{\mathcal{M}}
\newcommand{\Pc}{\mathcal{P}}
\newcommand{\Sc}{\mathcal{S}}

\newcommand{\cA}{\mathcal{A}}
\newcommand{\cC}{\mathcal{C}}
\newcommand{\cD}{\mathcal{D}}
\newcommand{\cF}{\mathcal{F}}
\newcommand{\cH}{\mathcal{H}}
\newcommand{\cI}{\mathcal{I}}
\newcommand{\cL}{\mathcal{L}}
\newcommand{\cM}{\mathcal{M}}
\newcommand{\cP}{\mathcal{P}}
\newcommand{\cQ}{\mathcal{Q}}

\newcommand{\ti}{{t_i}}
\newcommand{\tip}{{t_{i+1}}}
\newcommand{\rj}{{r_j}}
\newcommand{\rjp}{{r_{j+1}}}
\newcommand{\rjm}{{r_{j-1}}}

\newcommand{\D}{\mathbb{D}^{1,2}}
\newcommand{\ud}{\mathrm{d}}

\newcommand{\tY}{{\widetilde Y}}

\newcommand{\dtY}{\delta\widetilde Y}
\newcommand{\Xp}{X^{\pi}}
\newcommand{\dY}{{\delta Y}}
\newcommand{\dZ}{{\delta Z}}

\newcommand{\Yp}{Y^{\Re,\pi}}
\newcommand{\Zp}{Z^{\Re,\pi}}
\newcommand{\tYp}{{\widetilde Y^{\Re,\pi}}}
\newcommand{\bZp}{{\bar Z^{\Re,\pi}}}

\newcommand{\Lb}{\mathbf{L}}
\newcommand{\dtK}{\delta\widetilde{K}}
\newcommand{\dC}{{\delta C}}
\newcommand{\dF}{{\delta F}}

\newcommand{\alfam}{{{\check{\alpha}}}}
\newcommand{\tetam}{{{\check{\theta}}}}
\newcommand{\am}{{{\check{a}}}}
\newcommand{\Nm}{{{N^\am}}}
\newcommand{\dxi}{{\delta\xi}}

\makeatother

\begin{document}
\begin{frontmatter}

\title{Discrete-time approximation of multidimensional BSDEs with oblique reflections} %v24
\runtitle{Approximation of multidimensional reflected BSDEs}

\begin{aug}
\author[A]{\fnms{Jean-Francois} \snm{Chassagneux}\thanksref{t1}\ead[label=e1]{jean-francois.chassagneux@univ-evry.fr}},
\author[B]{\fnms{Romuald} \snm{Elie}\corref{}\thanksref{t2}\ead[label=e2]{elie@ceremade.dauphine.fr}}
and~%
\author[B]{\fnms{Idris}~\snm{Kharroubi}\ead[label=e3]{kharroubi@ceremade.dauphine.fr}}
\runauthor{J.-F. Chassagneux, R. Elie and I. Kharroubi}
\affiliation{Universit{\'e} d'Evry Val d'Essone and CREST, Universit{\'
e} Paris-Dauphine and CREST, and Universit{\'e} Paris-Dauphine and CREST}
\address[A]{J.-F. Chassagneux\\
D\'epartement de Math\'ematiques\\
Universit{\'e} d'Evry Val d'Essone\\
and CREST\\
France\\
\printead{e1}}
\address[B]{R. Elie\\
I. Kharroubi\\
CEREMADE, CNRS, UMR 7534\\
Universit{\'e} Paris-Dauphine\\
and CREST\\
France\\
\printead{e2}\\
\hphantom{E-mail: }\printead*{e3}}
\end{aug}

\thankstext{t1}{Supported by the ``Chaire Risque de cr\'edit,'' F\'ed\'eration Bancaire
Fran\c caise.}

\thankstext{t2}{Supported by the the Finance and Sustainable
Development Chair sponsored by EDF, CACIB and CDC.}

% HISTORY:
\received{\smonth{4} \syear{2010}}
\revised{\smonth{2} \syear{2011}}

% ABSTRACT
%
\begin{abstract}
In this paper, we study the discrete-time approximation of
multidimensional reflected BSDEs of the type of those presented by
Hu and Tang [\textit{Probab. Theory Related Fields}
\textbf{147} (2010) 89--121]
and generalized by Hamad\`ene and Zhang [\textit{Stochastic Process.
Appl.} \textbf{120} (2010) 403--426]. In comparison
to the penalizing approach followed by Hamad\`ene and Jeanblanc
[\textit{Math. Oper. Res.} \textbf{32} (2007) 182--192] or
Elie and Kharroubi [\textit{Statist. Probab. Lett.}
\textbf{80} (2010) 1388--1396], we study a more natural scheme based on
oblique projections. We provide a control on the error of the algorithm
by introducing and studying the notion of multidimensional discretely
reflected BSDE. In the particular case where the driver does not depend
on the variable~$Z$, the error on the grid points is of order
$\frac{1}{2} - \varepsilon$, $\varepsilon>0$.
\end{abstract}

% KEYWORDS
%

\end{frontmatter}

%s1 ###
\section{Introduction}

% Trying to optimize the productivity of a power station with start-up
%and switch-down costs,
%Hamadne and Jeanblanc~\cite{hamjea07} study the Snell envelope
%representation for optimal
%switching problems in continuous time. Observing that the difference
%between the
%two value functions starting in both modes of production is solution
%of a doubly
%reflected BSDE, they derive existence and uniqueness of solution to
%this problem.
%Nevertheless their approach restricts to optimal switching problems
%with only two possible
%modes of production. The extension to optimal switching problem in
%high dimension is studied
%by Carmona and Ludkovski~\cite{carlud05}, Porchet, Touzi and Warin
%by Pham, Ly Vath and Zhou~\cite{phalyvzho09} for an infinite time
%horizon consideration.
%In these papers, the resolution of optimal switching problem relies
%mostly on its link with
%systems of variational inequalities.
% The rigorous derivation of a multidimensional BSDE representation for
%this type of problem
%is obtained by Djehiche, Hamadne and Popier~\cite{djehampop07} and Hu
%and Tang~\cite{hutan08}. In this latter paper, they introduce the more
%general notion of multi-dimensional reflected BSDE, with dynamics of
%the form
%%%%

The main motivation of this paper is the discrete-time approximation
of the following system of reflected backward stochastic
differential equations (BSDEs)
%
%e1.1 ###
%
\begin{equation}\label{eqBSDECORIntro}\qquad
\cases{\displaystyle
Y ^i_t = g^i(X_T)+\int_t^T f^i(X_s,Y _s^i, Z _s^i)\,\ud s-\int
_t^T Z
^i_s \,\ud W_s\cr
\hphantom{Y ^i_t = }{} + K ^i_T - K ^i_t, &\quad$0 \le t \le
T$,\vspace*{2pt}\cr
\displaystyle Y ^i_t \geq\max_{j\in\Ic} \{Y ^j_t - c^{ij}(X_t)\}, &\quad
$0\leq
t\leq
T$,\vspace*{2pt}\cr
\displaystyle\int_0^T\Bigl[Y _t^i-\max_{j\in\Ic}\{Y ^j_t -
c^{ij}(X_t)\}\Bigr]
\,\ud K ^i_t=0,
&\quad$i \in
\cI$,}
\end{equation}
where $\Ic:=\{1,\ldots,d\}$, $f$, $g$ and $(c^{ij})_{i,j\in\Ic}$ are
Lipschitz functions and $X$ is the solution of a forward stochastic
differential equation (SDE).

These equations are linked to the solutions of \textit{optimal switching
problems}, arising, for example, in real option pricing. In the
particular case where $f$ does not depend on $(Y ,Z )$, a first study
of these equations was made by Hamad{\`e}ne and Jeanblanc \cite
{hamjea07}. They derive existence and uniqueness of solution to this
problem in dimension 2. The extension of this result to optimal
switching problems in higher dimension is studied by Djehiche, Hamad\`
ene and Popier~\cite{djehampop07}, Carmona and Ludkovski \cite
{carlud05}, Porchet, Touzi and Warin~\cite{portouwar09} or Pham, Ly
Vath and Zhou~\cite{phalyvzho09} for an infinite time horizon
consideration. In this last paper, the resolution of optimal switching
problems relies mostly on their link with systems of variational
inequalities. %\idris{on a l'impression que y a des bsde dans toutes
%ces ref}

Considering deterministic costs, Hu and Tang~\cite{hutan08} derive
existence and uniqueness of solution to this type of BSDE and relate it
to optimal switching problems between one-dimensional BSDEs. Extensions
developed in~\cite{hz08} and~\cite{chaelikha10} cover, in particular,
the existence of a unique solution to the BSDE~(\ref{eqBSDECORIntro}).
Recently two of the
authors related in~\cite{ek08} the solution
of~(\ref{eqBSDECORIntro}) to corresponding constrained BSDEs with
jumps. As presented in~\cite{ek09}, this type of BSDE can be
numerically approximated combining a penalization procedure with the
use of the backward scheme for BSDEs with jumps; see~\cite{boueli08}. Unfortunately, no
convergence rate is available for this algorithm. We present here a
more natural discretization scheme based on a geometric approach. For
any $t\le T$, all the components of the $Y _t$ process are
interconnected, so that the vector $Y _t$ lies in a random closed
convex set $\cQ(X_t)$ characterized by the cost functions
$(c^{ij})_{i,j\in\Ic}$. The vector process $Y $ is thus obliquely
reflected on the boundaries of the domain $\cQ(X)$ and we approximate
these continuous reflections numerically.

As in~\cite{majzha05,boucha08,cha08},
we first introduce a discretely reflected version of (\ref
{eqBSDECORIntro}), where the reflection occurs only on a deterministic
grid $\Re
=\{r_{0} := 0,\ldots,\break r_{\kappa}:=T\}\dvtx
Y^{\Re}_T = \tY^{\Re}_T := g(X_{T}) \in\cQ(X_{T})$, and, for $j
\le
\kappa-1$ and
$t \in[\rj,\rjp)$,
%
%e1.2 ###
%
\begin{equation} \label{BSDEDORintro}
\cases{
\displaystyle \tY^{\Re}_{t } = Y^{\Re}_\rjp+ \int_{t}^\rjp f(X_{u},
\tY^{\Re}_u, Z^{\Re}_u)\,\ud u - \int_{t}^\rjp Z^{\Re}_u \,\ud W_u,
\vspace*{2pt}\cr
\displaystyle Y^{\Re}_{t}= \tY^{\Re}_t\mathbf{1}_{ \{ t \notin\Re\} } +
\cP(X_t,\tY^{\Re}_t)\mathbf{1}_{ \{ t \in\Re\} } ,}
\end{equation}
where $\cP(X_t,\cdot)$ is the oblique projection operator on $\cQ(X_t)$,
for $t\le T$.
Extending the approach of Hu and Tang~\cite{hutan08}, we observe
that the solution to~(\ref{BSDEDORintro}) interprets as the value
process of a one-dimensional optimal BSDE switching problem with
switching times belonging to $\Re$. This allows us to prove a key
stability result for this equation.
We control the distance between $(Y^{\Re},Z^{\Re})$ and $(Y ,Z )$ in
terms of the mesh of the reflection grid. Due to the obliqueness of the
reflections, the direct argumentation of~\cite{boucha08,cha08} does
not apply. Using the reinterpretation in terms of switching BSDEs,
%as well as an innovative approach based on the introduction of a
%common dominating process,
we first prove that $Y^{\Re}$ approaches $Y $ on the grid points with
a convergence rate of order $\frac12-\eps$, $\eps>0$ uniformly in
$\Re
$, whenever the cost function is Lipschitz and~$f$ is bounded in $z$
(see Theorem~\ref{ThmErrorDiscretize}). Imposing more regularity on
the cost functions, we control the convergence rate of $(Y^{\Re
}_t,Z^{\Re}_t)_{0\le t \le T}$ to $(Y _t,Z _t)_{0\le t \le T}$ (see
Theorem~\ref{ThmErrorDiscretize2}). % when the costs are smoother
%(say $C^{2}_{b}$), see Theorem~\ref{ThmErrorDiscretize2}.

%Whenever $f$ is bounded in the variable $Z$,
%this particular feature leads to a convergence rate of order $|
%obliquely reflected BSDE \reff{BSDE DORintro} and the continuously
%reflected one \reff{eq BSDE COR Intro}. This convergence rate is of
%order $|\Re|^{{1 \over2}-\eps}$ when we restrict to the grid points.

We then consider a Euler type approximation scheme associated to the
BSDE~(\ref{BSDEDORintro}) defined on $\pi=\{t_0,\ldots,t_n\}$ by
$\Yp_T:=g(\Xp_T)$ and, for $i\in\{n-1,\ldots,0\}$,
%
%e1.3 ###
%
\begin{equation}\label{schemeintro}
\cases{
\bZp_{\ti} := (\tip-\ti)^{-1}\mathbb{E}[\Yp_\tip(W_\tip
-W_\ti)'\mid\mathcal{F}_{\ti}] ,
\vspace*{2pt}\cr
\tYp_{\ti} := \mathbb{E}[\Yp_\tip\mid\mathcal{F}_{\ti
}] + (\tip- \ti)
f(\Xp_\ti,\tYp_\ti,\bZp_{\ti}) ,
\vspace*{2pt}\cr
\Yp_{\ti} := \tYp_\ti\mathbf{1}_{ \{ \ti\notin\Re\} } + \cP
(\Xp_\ti,
\tYp_\ti)\mathbf{1}_{ \{ \ti\in\Re\} } ,}
\end{equation}
where $\Xp$ is the Euler scheme associated to $X$. It is now well known
(see, e.g.,~\cite{boutou04,zha04}), that the convergence rate of the
scheme~(\ref{schemeintro}) to the
solution of~(\ref{BSDEDORintro}) is controlled by the regularity of
$(Y,Z)$ through the quantities
\[
\mathbb{E}\biggl[\sum_{i <n} \int_{\ti}^{\tip} |Y^{\Re}_{t} - Y^{\Re
}_{\ti}|^2 \,\ud t\biggr] \quad\mbox{and}\quad \mathbb{E}\biggl[\sum_{i <n}
\int_{\ti}^{\tip} |Z^{\Re}_{t} - \bar Z^{\Re}_{\ti}|^2 \,\ud t\biggr]
\]
with $\bar Z^{\Re}_{\ti} = \frac{1}{\tip-\ti}\mathbb{E}
[\int_{\ti}^\tip Z^{\Re}_{t} \,\ud t\mid\mathcal{F}_{\ti}
]$, for $i\le n$.

Using classical Malliavin differentiation tools, we prove a
representation for~$Z^{\Re}$, extending the results of \cite
{boucha08,cha08} to the system of discretely reflected BSDEs (\ref
{BSDEDORintro}). We deduce the expected regularity results on $(Y^{\Re
},Z^{\Re})$ and, using the techniques of~\cite{cha09}, Chapter 3, we
obtain in a very general setting the convergence of~(\ref{schemeintro})
to~(\ref{BSDEDORintro}).
However, due to the obliqueness of the reflections, the projection
operator $\Pc(X,\cdot)$ is only $L_\cP$-Lipschitz with $L_\cP:=\sqrt{d}>1$,
leading to a convergence rate controlled by $|L_\cP|^{\kappa}( |\pi
|^{1/4} + \kappa^{1/2}|\pi|^{1/2})$, where we recall that
$\kappa$ is the number of points in the reflection grid $\Re$. The term
$|L_\cP|^{\kappa}$ can be very large even for small $\kappa$ and leads
to a poor logarithmic convergence rate when passing to the limit
$\kappa\rightarrow\infty$ for the approximation of~(\ref
{eqBSDECORIntro}). In the particular case where $f$ does not depend on $z$,
we are
able to get rid of the $|L_\cP|^\kappa$ term. %Our innovative
%approach relies on the reinterpretation of the solution to
%combined with
%the introduction of a convenient auxiliary process dominating both
%solutions of \reff{BSDE DORintro} and \reff{scheme intro}.

Our innovative approach relies on the use of comparison results to get
a~control of the involved quantities:
\begin{itemize}
\item we interpret the solution of~(\ref{BSDEDORintro}) as a value
process of an optimization problem, which allows us to get a control of
the distance between the continuously and discretely reflected BSDEs;

\item we introduce a convenient auxiliary process dominating both
solutions~(\ref{BSDEDORintro}) and~(\ref{schemeintro}) to get a control
of the distance between these quantities.
\end{itemize}
Combining the previous estimates, we deduce
the convergence of the discrete time scheme~(\ref{schemeintro}) to
the solution of~(\ref{eqBSDECORIntro}) with a convergence rate of
order $\frac12-\eps$, $\eps>0$, on the grid\vadjust{\goodbreak} points, whenever $\Re
=\pi$
and $f$ is independent of $Z$. Whenever the cost functions are
constant, all the previous estimates hold true with $\eps=0$. We want
to emphasize that all these results are obtained without any assumption
on the nondegeneracy of the volatility matrix $\sigma$.

The rest of the paper is organized as follows. In Section~\ref{sedisc},
we introduce the notion of discretely obliquely reflected BSDEs,
connect it with optimal switching problems and give the fundamental
stability result. Section~\ref{secreg} focuses on the regularity of the
solution to this new type of BSDE. This analysis leads to precious
estimates allowing us to deduce the convergence of the associated
discrete time scheme (see Section~\ref{SubSecScheme}). Afterward,
Section~\ref{secontinuous} focuses on the extension to the continuously
reflected case and provides a~convergence rate of the discretely
reflected BSDE to the continuously one, whenever the driver $f$ is
bounded in the variable~$Z$. The global error of the scheme is provided
at the end of this section. Some a priori estimates are reported in
the \hyperref[app]{Appendix}.

% Section 3 focuses on the discrete time approximation of discretely
%reflected BSDEs. We present the discrete time scheme and study its
%convergence. In the particular case where the driver function does not
%involve the variable $Z$, we provide a rate of convergence. Finally,
%Section 4 study the extension to the continuously reflected case. In
%the case where the driver $f$ is bounded in the variable $Z$, we
%provide a rate of convergence of the discretely reflected BSDE to the
%continuously one.

\subsection*{Notation} Throughout this paper we are given a finite time
horizon $T$ and a probability space $(\Omega,\Fc,\PP)$ endowed with a
$d$-dimensional standard Brownian motion $W =(W_t)_{t\geq0}$. The
filtration $\F=(\Fc_t)_{t\le T}$ generated by the Brownian motion is
supposed to satisfy the usual conditions. Here, $\mathfrak{P}$ denotes
the $\sigma$-algebra on $[0,T]\times\Omega$ generated by $\F
$-progressively measurable processes. Any element $x \in\R^\ell$ with
$\ell\in\N$ will be identified to a column vector with $i$th component
$x^i$ and Euclidean norm $|x|$. For $x,y\in\R^\ell$, $x\cdot y$ denotes the
scalar product of $x$ and~$y$, and $x'$ denotes the transpose of $x$.
We denote by $\succeq$ the component by component partial ordering
relation on vectors.
$\Mc^{m,d}$ denotes the set of real matrices with $m$ lines and $d$ columns.
We denote by $C^k_{b}$ the set of functions from $\R^d$ to $\R$ with
continuous and bounded derivatives up to order $k$. For a function
$f\in C^1$, $\nabla_x f$ denotes the Jacobian matrix of $f$ with
respect to $x$. For ease of notation, we will sometimes write $\E
_t[\cdot]$ instead of $\E[\cdot| \Fc_t]$, $t\in[0,T]$. In the
following, we shall use the notation without specifying the dimension
nor the dependence in $\omega\in\Omega$ when it is clearly given by
the context.
Finally, for any $p\ge1$, we introduce the following:

\begin{itemize}
\item
the set ${\Sc^p}$ of real-valued \textit{c\`adl\`
ag}\setcounter{footnote}{1}\footnote
{French acronym meaning right continuous with left limit.}
%gauche} meaning right continuous with left limits.} (resp. continuous)
$\mathfrak{P}$-measurable processes $Y=(Y_t)_{0\leq t\leq T}$
satisfying
$\|Y\|_{_{{\Sc^p}}} := \mathbb{E}[\sup_{0\leq t\leq T} |Y_t|^p]
^{1/p}
< \infty$.

\item the set ${\Hc^p}$ of $\R^d$-valued $\mathfrak{P}$-measurable
processes $Z=(Z_t)_{0\leq t\leq T}$ such that
$\|Z\|_{_{\Hc^p}} := \mathbb{E}[(\int_0^T |Z_t|^2 \,\ud t)^{p/2}
]^{1/p} < \infty$.

\item the closed subset $\mathbf{A}^p$ of ${\Sc^p}$ consisting of
nondecreasing processes $K$ satisfying $K_0=0$.

In the sequel we denote by $C_L$ a constant whose value may change
from line to line
but which depends only on $L$. We use the notation $C^p_L$ whenever
it depends on some other parameter $p>0$.\vadjust{\goodbreak}
\end{itemize}

%A mettre dans notation:(here $\mathfrak{P}$ denote the $
%measurable processes)

%%%%%%%%%%%%%%%%%%%%%%%%%%
%%%%%%%%%%%%%%%%%%%%%%%%%%
%%%%%%%%%%%%%%%%%%%%%%%%%%
%%%%%%%%%%%%%%%%%%%%%%%%%%

%!TEX root = main.tex

%s2 ###
\section{Discretely obliquely reflected BSDE}
\label{sedisc} In the beginning of this section we define and
study discretely obliquely reflected BSDEs in a general setting. In
particular, we show how their solutions relate to the solutions of
one-dimensional optimal switching problems, where the switching
times are restricted to lie in a discrete time set. This allows us to
prove a stability result for obliquely RBSDEs which will be used
several times in the paper.
%Finally, at the end of the section, we
%introduce the diffusion framework we are going to use in the sequel.
%In the next section, we derive regularity results for discretely
%reflected BSDEs within this diffusion framework.

%This formulation allows to characterize its proximity with the
%corresponding continuously obliquely reflected BSDE. We prove that the
%${\bf L^2}$ distance between both solution is of the order of the
%square root of the mesh size of the discretization grid. This leads to
%the presentation of a convergent discrete-time scheme for
%multidimensional reflected BSDE.

%s2.1 ###
\subsection{Definition}

A discretely obliquely reflected BSDE is a reflected BSDE
%of type \eqref{eqBSDECORIntro}
where the reflection is only allowed on a
discrete time set.

We thus consider a grid $\Re:=\{ r_{0}=0,\ldots, r_{\kappa}=T \}$ of
the time interval $[0,T]$ satisfying
%
%e2.1 ###
%
\begin{equation}
\label{CondRe}
|\Re|:=\max_{1\leq k\leq\kappa}|r_{k}-r_{k-1}|\le\frac{L}{\kappa}.
\end{equation}

%Let us precise the particular constraint imposed to the value of the
%BSDE.

We also consider a matrix valued process $C = (C^{ij})_{1\leq i,j\leq
m}$ %\in\HP{2}(M^{d})$
such that $C^{ij}$ belongs to $\Sc^2$ for $i,j\in\{1,\ldots,d\}$ and
satisfies the structure condition
%
%e2.2 ###
%
\begin{equation}\label{CostStructureabs}\qquad
\cases{
C^{ii}_{t}=0 , &\quad for $1\le i \le d$ and $0\leq t\leq T$;\vspace*{2pt}\cr
\displaystyle \inf_{0\leq t\leq T} C^{ij}_{t}> \frac1L ,
&\quad for $1\le i,j \le d$ with $i\neq j$; \vspace*{2pt}\cr
\displaystyle \inf_{0\leq t\leq T}C^{ij}_{t} + C^{jl}_{t} - C^{il}_{t}>0,
&\quad for $1\le i,j,l\le d$ with $i\neq j, j\neq l$.}
\end{equation}
%
% \jef{regularite integrabilite des cout pour la representation a
%verifier}
%d'etre uniforme en $\omega$?}
We introduce a random closed convex set family associated to $C$:
\[
\cQ_{t} := \Bigl\{ y \in\R^d \mid y^i \geq\max_{j} (y^j -
C_{t}^{ij} ) , 1\le i \le d \Bigr\},\qquad 0\leq t\leq T,
\]
and the oblique projection operator onto ${\cQ_{t}}$,
denoted $\cP_{t}$ and defined by
\[
\cP_{t} \dvtx y \in\R^d \mapsto\Bigl( \max_{j\in\Ic} \{
y^j -
C^{ij}_{t} \} \Bigr)_{1\le i\le d},
\]
which is $\mathfrak{P}\otimes\mathcal{B}(\R^d)$-measurable.
\begin{Remark}\label{repropP}
(i) It follows from the structure condition~(\ref{CostStructureabs})
that $\cP$ is increasing with respect to the partial ordering
relation $\succeq$, where $y\succeq y'$ means $y^i\geq(y')^i$ for
all $i\in\cI$.

(ii) An easy calculation leads to
\[
| \cP_t(y_1) - \cP_t(y_2) | \le\sqrt{d} |y_1 - y_2|\qquad
\mbox{for any } y_1,y_2 \in\R^d.
\]
We observe that the constant $\sqrt{d}$ is optimal in our setting
taking, for example, $y_{1} := (\max_{i,j} C_t^{ij},0,\ldots,0)$ and
$y_{2} := (\max_{i,j} C_t^{ij}+1,0,\ldots,0)$.
Thus $\cP_t$ is $L_\cP$-Lipschitz continuous with $L_\cP:= \sqrt{d}$.
%, and, due to the obliqueness of the reflexions, one easily checks
%that $L_\cP>1$. Hence, we recall for later use that
%| \cP_t(y_1) - \cP_t(y_2) | & \le L_\cP|y_1 - y_2| ,
% \mbox{with } 1< L_\cP\le\sqrt{d} .
\end{Remark}

Finally, we are also given a random variable $\xi\in[L^{2}(\cF
_{T})]^d$ valued in $\cQ_{T}$, representing the
terminal value of the BSDE and a random function
$F\dvtx\Omega\times[0,T]\times\R^d\times\cM^{d,q}\rightarrow\R^d$ which
is
$\mathfrak{P}\otimes\mathcal{B}(\R^d)\otimes\mathcal{B}(\cM
^{d,q})$-measurable\vadjust{\goodbreak}
and satisfies the Lipschitz property
\[
|F(t,y,z)-F(t,y',z')| \le L(|y-y'|+|z-z'|)
\]
for all $(t,y,y',z,z')\in[0,T]\times(\R^d)^2\times(\cM^{d,q})^2$,
$\PP$-a.s.
We shall also assume that

\begin{longlist}[(HF)]
\item[(HF)] The component $i$ of $F(t,y,z)$ depends only on the
component $i$ of the vector $y$ and on the row $i$ of the matrix $z$,
that is, $F^i(t,y,z)=F^i(t,y^i,z^i)$.
\end{longlist}

Given\vspace*{1pt} this set of data $(\Re,C,F,\xi)$, a discretely obliquely
reflected BSDE, denoted $\cD(\Re,C,F,\xi)$, is a triplet
$(\tY^{\Re},Y^{\Re},Z^{\Re})\in({\Sc^2}\times{\Sc^2}\times\Hc
^2)^\Ic$
satisfying $Y^{\Re}_T
= \tY^{\Re}_T := \xi\in\cQ_{T}$, and defined in a backward manner,
for $j \le\kappa-1$ and $t \in[\rj,\rjp)$, by
%
%e2.3 ###
%
\begin{equation}\label{eqdeBSDEDOR1abs}
\cases{
\displaystyle \tY^{\Re}_{t } = Y^{\Re}_\rjp+ \int_{t}^\rjp F(u, \tY^{\Re}_u,
Z^{\Re}_u)\,\ud u - \int_{t}^\rjp Z^{\Re}_u \,\ud W_u , \cr
\displaystyle Y^{\Re}_{t}= \tY^{\Re}_t\mathbf{1}_{ \{ t \notin\Re\} } +
\cP_{t}(\tY^{\Re}_t)\mathbf{1}_{ \{ t \in\Re\} } .}
\end{equation}
This rewrites equivalently for $t \in[0,T]$ as
%
%e2.4 ###
%
\begin{equation} \label{eqBSDEDORabs}
\cases{
\displaystyle \tY^{\Re}_{t } =
\xi+ \int_{t}^T F(u,\tY^{\Re}_u,Z^{\Re}_u)\,\ud u -
\int_{t}^T Z^{\Re}_u \,\ud W_u + (K ^{\Re}_T-K ^{\Re}_t) , %
% 0
\vspace*{2pt}\cr
\displaystyle K ^{\Re}_t := \sum_{r \in\Re\setminus\{ 0 \}} \Delta K ^{\Re
}_r
\mathbf{1}_{ \{ r \le t \} } , \cr
\qquad\mbox{with $\Delta K ^{\Re}_t := Y ^{\Re
}_t - \tY
^{\Re}_t =
-(\tY^{\Re}_t -
\tY^{\Re}_{t-})$}.}
% 0\le t \le T.
\end{equation}
Observe that $K ^{\Re}\in(\mathbf{A}^2)^\Ic$, since $C^{ij}$ is
nonnegative and valued in $\Sc^2$, for any $i,j\in\Ic$.

%K _t := \sum_{r \in\Re\setminus\set{0}}
%-(\tY_t - \tY_{t-}) , 0\le t \le T.

We shall also use the following integrability condition
for some $p\ge2$:
{\renewcommand{\theequation}{C$_p$}
\begin{equation}\label{Cp}
|\xi|^p + \sup_{t\in[0,T]}
|C_{t}|^p + \int_0^T |F(s,0,0)|^p \,\ud s \le\beta,
\end{equation}}

\vspace*{-4pt}

\noindent where $\beta$ is a positive random variable satisfying $\mathbb
{E}[\beta]\le
C_{L}$. Importantly, $\beta$ does not depend on $\Re$.

The proof of the following a priori estimates is postponed until
the \hyperref[app]{Appendix}.

\begin{Proposition}\label{estprioedsrdrabs}
Assume that~(\ref{Cp}) holds for some given $p \ge2$, there
exists a unique solution $(\tY^{\Re},Y^{\Re},Z^{\Re})$ to (\ref
{eqdeBSDEDOR1abs}) and
it satisfies
\[
%C^p_L ,
\|\tY^{\Re}\|_{_{ \Sc^p}}+\|Z^{\Re}\|_{_{ \Hc^p}}+\|K ^{\Re}_T\|
_{_{\mathbf{
L}^p}} \le C^p_L.
\]
%
%recalling that $C^p_L$ does
%not depend on $\Re$.
\end{Proposition}

%a standard BSDE, existence and uniqueness follow from a
%concatenation procedure and~\cite{parpen90}.

%The integrability result is a direct application of Proposition

%s2.2 ###
\subsection{Corresponding optimal switching problem}
\label{subsediscopt}

In this subsection, we interpret the solution of the discretely
obliquely RBSDE~(\ref{eqBSDEDORabs}) as the value process of a
corresponding optimal switching problem, where the possible
switching times are restricted to belong to the grid $\Re$. Our
approach relies on similar arguments as the one followed by Hu and
Tang~\cite{hutan08} in a framework with continuous reflections.\vadjust{\goodbreak}

A switching strategy $a$ is a nondecreasing sequence of stopping times
$(\theta_j)_{j\in\N}$, combined with a sequence of random variables
$(\alpha_j)_{j\in\N}$ valued in~$\Ic$, such that $\alpha_j$ is
$\Fc
_{\theta_j}$-measurable, for any $j\in\N$. %For such a strategy $a$
We denote by $ \cA$ the set of such strategies.
For $a= (\theta_{j},\alpha_{j})_{j\in\N}\in\Ac$, we introduce
$N^a$ the
(random) number of switches before $T$ as
%
%e2.5 ###
%
\setcounter{equation}{4}
\begin{equation} \label{eqdefNa}
N^a = \#\{k\in\N^*\dvtx\theta_{k} \leq T \}.
\end{equation}
To any switching strategy $a=(\theta_{j},\alpha_{j})_{j\in\N}\in
\cA$, we
associate the current state process $(a_t)_{t\in[0,T]}$ and the
compound cost process $(A^a_t)_{t\in[0,T]}$ defined, respectively, by
\[
a_t := \alpha_{0}{\mathbf1}_{\{0\leq t<\theta_0\}}+\sum
_{j=1}^{N^a} \alpha
_{j-1} {\mathbf1}_{\{\theta_{j-1}\leq t < \theta_{j}\}}
\quad\mbox{and}\quad A^a_t
:= \sum
_{j=1}^{N^a} C^{\alpha_{j-1}\alpha_{j}}_{\theta_{j}} {\mathbf1}_{\{
\theta_j\le
t \le T\}} % 0\le t\le T .
\]
for $0\le t\le T$.
%An admissible switching strategy $a$ is a switching strategy satisfying
%$
For $(t,i)\in[0,T]\times\cI$, the set $\cA_{t,i}$ of admissible
strategies starting from $i$ at time $t$ is defined by
\[
\cA_{t,i} = \{ a=(\theta_{j},\alpha_{j})_{j}\in\Ac\mid\theta_{0}=
t,\alpha_{0}= i,\mathbb{E}[|A_{T}^a|^2]< \infty\}.
\]
%
% a nondecreasing sequence of stopping times $(\theta_j)_{j\in\N}$
%valued in $\Re$, combined with a sequence of random variables $(
%is $\Fc_{\theta_j}$-measurable, for any $j\in\N$.
%To any switching strategy $a=(\theta_{j},\alpha_{j})_{j\in\N}$, we
%associate the current state process $(a_t)_{t\in[0,T]}$ and the
%compound cost process $(A^a_t)_{t\le T}$ defined respectively by
% \[
% a_t := \sum_{j\in\N} \alpha_j {\mathbf1}_{\{\theta_j\le t <
% .
% \]
Similarly, we introduce $\Ac^\Re_{t,i}$, the restriction to $\Re
$-admissible strategies
\[
\Ac^\Re_{t,i} := \{ a=(\theta_j,\alpha_j)_{j\in\N} \in\cA
_{t,i} \mid \theta_{j}\in\Re, \forall j\leq N^a \}%
% 0\le t \le T , i \le d .
\]
and denote $\Ac^\Re:=\bigcup_{i\le d} \Ac^\Re_{0,i}$.
%and by $\Ac^\Re_{t,i}$ the set of admissible strategies starting from
%$i$ at time $t$ :
%is denoted and we define
% | (\theta_0,\alpha_0) = (t,i) \} , 0\le t \le T ,
% i \le d .

For $(t,i)\in[0,T]\times\cI$ and $a\in\Ac^\Re_{t,i}$, we
consider as
in~\cite{hutan08}
the associated one-dimensional switched BSDE defined by
%
%e2.6 ###
%
\begin{eqnarray}\label{eqUa}
U^a_u &=& \xi^{a_T} + \int_u^T F^{a_s}(s,U^a_s,V^a_s)\,\ud s - \int_u^T
V^a_s \,\ud W_s\nonumber\\[-8pt]\\[-8pt]
&&{} - A^a_T + A^a_u ,\qquad t\le u \le T.\nonumber
\end{eqnarray}

Theorem 3.1 in~\cite{hutan08} interprets each component of the
solution to the continuously reflected BSDE~(\ref{eqBSDECORIntro}) as
the Snell envelope associated to switched processes of the form (\ref
{eqUa}), where the switching strategies $a$ are not restricted to lie
in the reflection grid $\Re$.
The next theorem is a new version of this Snell envelope
representation adapted to the context of discretely obliquely reflected
BSDE~(\ref{eqBSDEDORabs}).
%
% \romu{Essayer de dtailler l'interpretation en terme de probleme de
%switch}
%
\begin{Theorem}\label{THswR}
Assume that \textup{(C$_2$)} is in force.
For any $i\in\cI$ %\{1,\ldots, d\}$
and $t\in[0,T]$, the following hold:
\begin{longlist}
\item The process $\tY^{\Re}$ dominates any $\Re$-switched BSDE,
that is,
%
%e2.7 ###
%
\begin{equation}\label{eqUpperboundYU}
U^a_t \le(\tY^{\Re}_t)^i ,\qquad\mathbb{P}\mbox{-a.s.}
\mbox{ for any }
a\in
\Ac^\Re_{i,t}.
\end{equation}
\item Define the strategy $a^*=(\theta_j^*,\alpha_j^*)_{j\ge0}$
recursively by $(\theta^*_0,\alpha^*_0):=(t,i)$ and, for $j\ge1$,
\begin{eqnarray*}
\theta^*_j &:=& \inf\Bigl\{ s\in[\theta^*_{j-1},T]\cap\Re
\bigm|
(\tY_s^{\Re})^{\alpha_{j-1}^*} \le\max_{k\neq\alpha^*_{j-1}}
\{ (\tY^\Re_s)^{k} - C^{\alpha^*_{j-1}k}_s \} \Bigr\},\\[-2pt]
\alpha^*_j &:=& \min\Bigl\{ \ell\neq\alpha^*_{j-1} \bigm|
(\tY^\Re_{\theta_j^*})^\ell- C^{\alpha^*_{j-1}\ell}_{\theta
^*_{j}} =
\max_{k\neq\alpha^*_{j-1}} \{ (\tY^{\Re}_s)^k -
C^{\alpha^*_{j-1} k}_{\theta_{j}^*} \} \Bigr\}.
\end{eqnarray*}
Then, we have $a^*\in\cA^\Re_{t,i}$ and
%
%e2.8 ###
%
\begin{equation}
(\tY^{\Re})^i_t = U^{a^*}_t ,\qquad\mathbb{P}\mbox{-a.s.}
\end{equation}

\item The following ``Snell envelope'' representation holds:
%
%e2.9 ###
%
\begin{equation}\label{SnellDiscrete}
(\tY^{\Re})^{i}_t = \mathop{\esssup}_{a\in\Ac^\Re_{t,i}} U^a_t,\qquad
\mathbb{P}\mbox{-a.s.}
\end{equation}
\end{longlist}
%
% We call $a^*$ the optimal strategy associated to $\tY$ in $\cA^
\end{Theorem}
\begin{pf}
Observe first that assertion (iii) is a direct
consequence of (i) and (ii). Let us fix $t\in[0,T]$ and
$i\in\cI$.\vspace*{8pt}

\textit{Step} 1. We first prove (i). %\{1,\ldots,d\}$.

Set $a=(\theta_k,\alpha_k)_{k\ge0} \in\Ac^\Re_{t,i}$ and the process
$(\tY^a,Z^a)$ defined, for $s\in[t,T]$, by
%
%e2.10 ###
%
\begin{eqnarray}\label{defYaZa}
\tY^a_s &:=& \sum_{k\ge0} (\tY_s^{\Re})^{\alpha_{k}} {\mathbf1}_{\{
\theta_k \le
s < \theta_{k+1}\}} + \xi^{a_T}{\mathbf1}_{\{s=T\}}\quad
\mbox{and}\nonumber\\[-9pt]\\[-9pt]
Z^a_s &:=& \sum_{k\ge0} (Z_s^{\Re})^{\alpha_{k}} {\mathbf1}_{\{
\theta_k \le s <
\theta_{k+1}\}} .\nonumber
% K _s := \sum_{k\ge0} \tY_s^{\alpha_{k}} {\mathbf1}_{\theta_k
\end{eqnarray}
Observe that these processes jump between the components of the
discretely reflected BSDE~(\ref{eqdeBSDEDOR1}) according to the
strategy $a$, and, between two jumps, we have
%
%e2.11 ###
%
\begin{eqnarray}\label{babel}\qquad
\tY^{a}_{\theta_{k}}
% =
% \tY^{\Re,\alpha_{k}}_{\theta_{k}}
&= & (Y^{\Re}_{\theta_{k+1}})^{\alpha_{k}} + \int
_{\theta_{k}}^{\theta_{k+1}} F^{\alpha_k}(s,(\tY_s^{\Re
})^{\alpha
_{k}},(Z^{\Re}_s)^{\alpha_k}) \,\ud s - \int_{\theta
_{k}}^{\theta
_{k+1}} (Z^{\Re}_s)^{\alpha_k}\,\ud W_s\nonumber\\[-2pt]
&&{} + (K ^{\Re
}_{\theta
_{k+1}-})^{\alpha_{k}} - (K ^{\Re}_{\theta_{k}})^{\alpha
_{k}}\nonumber\\[-9pt]\\[-9pt]
&= &
\tY^{a}_{\theta_{k+1}} + \int_{\theta_{k}}^{\theta_{k+1}}
F^{a_s}(s,\tY
^a_s,Z^a_s) \,\ud s - \int_{\theta_{k}}^{\theta_{k+1}} Z^a_s \,\ud W_s
+ (K ^{\Re}_{\theta_{k+1}-})^{\alpha_{k}}\nonumber\\[-2pt]
&&{} - (K ^{\Re}_{\theta
_{k}})^{\alpha_{k}}
% + (K ^{\alpha_{k+1}}_{\theta_{k+1}} - K ^{\alpha_{k+1}}_{
+ \bigl((Y^{\Re}_{\theta_{k+1}})^{\alpha_{k}} - (\tY^{\Re}_{\theta
_{k+1}})^{\alpha_{k+1}}\bigr) ,\qquad k\ge0.\nonumber
\end{eqnarray}
Introducing
\begin{eqnarray*}
K ^a_s&:=&\sum_{k = 0}^{N^{a}-1} \biggl[\int_{(\theta_{k}\wedge s,
\theta
_{k+1}\wedge s)} \mathrm{d}(K ^{\Re}_u)^{\alpha_k}\\[-2pt]
&&\hspace*{25.8pt}{}+\mathbf{1}_{\{ \theta
_{k+1}\leq s \}}\bigl((Y^{\Re}_{\theta_{k+1}})^{\alpha_k}-(\tY
^{\Re
}_{\theta_{k+1}})^{\alpha_{k+1}}+C^{\alpha_{k}\alpha_{k+1}}_{\theta
_{k+1}}\bigr)\biggr]
\end{eqnarray*}
for $s\in[t,T]$, and summing up~(\ref{babel}) over $k$, we get, for $
t\leq u \leq
T$,
\[
% \[
\tY^{a}_u
=
\xi^{a_T} + \int_{u}^{T} F^{a_s}(s,\tY^a_s,Z^a_s) \,\ud s - \int_{u}^{T}
Z^a_s \,\ud W_s
- A^a_T+A^a_u + K ^a_T -K ^a_u.
\]
Using the relation $Y^{\Re}_{\theta_{k}}=\cP_{\theta_{k}}(\tY
^{\Re}_{\theta_{k}})$ for all $k\in\{0,\ldots,N^{a}\}$, we check that
$K ^a$ is increasing.
% Since $K ^a$ is increasing, $Y_{t}\in\cQ(X_{t})$ for all $t\in\Re$
%and
Since $U^a$ solves~(\ref{eqUa}), we deduce by a comparison argument
(see~\cite{pen99}, Theorem 1.3) that $U^a_{t} \le\tY^a_{t}$.
Since $a$ is arbitrary in $\Ac^\Re_{t,i}$, we deduce~(\ref{eqUpperboundYU}).\vspace*{8pt}

\textit{Step} 2. We now prove (ii).

Consider the strategy $a^*$
given above as well as the associated process $(\tY^{a^*},Z^{a^*})$
defined as in~(\ref{defYaZa}). By definition of $a^*$, we have
\[
(Y_{\theta^*_{k+1}}^{\Re})^{\alpha^*_{k}} = (\Pc_{\theta
^*_{k+1}}(\tY^{\Re}_{\theta^*_{k+1}}))^{\alpha^*_{k}} =
(\tY
^{\Re}_{\theta^*_{k+1}})^{\alpha^*_{k+1}} - C^{\alpha^*_{k}\alpha
^*_{k+1}}_{\theta^*_{k+1}} ,\qquad k\ge0,
\]
which gives
%
%e2.12 ###
%
\begin{equation}
\int_{(\theta^*_{k}, \theta^*_{k+1})}
\mathrm{d}(K ^{\Re}_s)^{\alpha^*_{k}}=0 \quad\mbox{and}\quad
(Y^{\Re}_{\theta^*_{k+1}})^{\alpha^*_{k}}-(\tY^{\Re}_{\theta
^*_{k+1}})^{\alpha^*_{k}}+C^{\alpha^*_{k}\alpha^*_{k+1}}_{\theta^*_{k+1}}=0\hspace*{-35pt}
\end{equation}
for all $k\in\{0,\ldots,N^{a^*}-1\}$.
% Since $K $ is constant between the jump times of $a^*$, we deduce
We deduce from~(\ref{CostStructureabs}) that
\[
\tY^{a^*}_u
=
\xi^{a^*_T} + \int_{u}^{T} F^{a^*_s}(s,\tY^{a^*}_s,Z^{a^*}_s) \,\ud s -
\int_{u}^{T} Z^{a^*}_s \,\ud W_s
- A^{a^*}_T+A^{a^*}_u ,\qquad t\leq u\leq T.
\]
Hence, $(\tY^{a^*}, Z^{a^*})$ and $(U^{a^*}, V^{a^*})$ are solutions
of the same BSDE and $(\tY^{\Re}_t)^i = U^{a^*}_t$.
To complete the proof, we only need to check that $a^*\in\Ac^\Re$,
that is, $\E|A^{a^*}_{T}|^2<\infty$.
By definition of $a^*$ on $[t,T]$ and the structure condition on the
cost~(\ref{CostStructureabs}), we have $|A^{a^*}_t|\leq\max_{k\neq
i}|C^{i,k}_{t}|$ which gives $\E[|A^{a^*}_t|^2] \leq C_{L}$.
Combining
\[
A^{a^*}_T = \tY^{a^*}_T-\tY^{a^*}_t + \int_{t}^{T} F^{a^*_s}(s,\tY
^{a^*}_s,Z^{a^*}_s) \,\ud s-\int_{t}^{T} Z^{a^*}_s \,\ud W_s+ A^{a^*}_t
\]
with the\vspace*{1pt} Lipschitz property of $F$ and the fact that $(\tY^{\Re
},Z^{\Re
})\in({\Sc}^2\times{\Hc}^2)^\Ic$ (recall Proposition \ref
{estprioedsrdrabs}), we get the square integrability of $A_{T}^{a^*}$
and the proof is complete.
\end{pf}
%
%{\it(iii)} Since $a^*\in\Ac^\Re_{t,i}$, estimate \reff{Snell
%Discrete} is a direct consequence of (i) and (ii).
%continuite par
%rapport aux parametres pour pouvoir les regulariser. on pourra aussi
%commenter
%sur le 'compare to control'}
%le cas ou f n'est pas markovien
%
\begin{Remark}\label{reTHswRone}
Although the optimal strategy $a^*$ depends on the initial
parameters $t$ and $i$, we omit the script $(t,i)$ for ease of
notation. %since its proof does not use the Markovian property of the
%BSDE. This remark will be used in the sequel.
\end{Remark}
%
%{\rm Notice that Theorem~\ref{THswR} holds true for a random generator
%$(F(t,\cdot,\cdot))_{t}$ such that $\esp{\int_{0}^T|F(s,0,0)|^2ds}$ $<$ $
%%since its proof does not use the Markovian property of the BSDE. This
%remark will be used in the sequel.

Combining the previous representation with the a priori estimates of
Proposition~\ref{estprioedsrdrabs} and the structure condition
(\ref{CostStructureabs}), we deduce the following estimates, whose
proof is postponed until the \hyperref[app]{Appendix}.
\begin{Proposition}\label{corhatA} Assume that~(\ref{Cp})
holds for some given $p \ge2$, then
%%\esp{\sup_{t\in[0,T]}|U^{a^*}_{t}|^p+\int_{0}^T|V^{a^*}_{s}|^p\,\ud s
%+
%|A^{a^*}_{T}|^2 } \le C^p_L ,
%%\|U^{a^*}\|_{_{{\bf\Sc^p}}}+\|V^{a^*}\|_{_{{\bf\Hc^p}}} + \|A^{a^*}
%+ \E[|N^{a^*}|^p] \le C^p_L ,p \geq2 ,%(t,i)\in[0,T]\times\Ic
%
\[
\mathbb{E}\biggl[\sup_{s\in[t,T]} |U^{a^*}_{s}|^p+ \biggl(\int
_{t}^T|V^{a^*}_{u}|^2\,\ud u\biggr)^{p/2} + |A^{a^*}_{T}|^p +
|N^{a^*}|^p \biggr] \le C^p_L
\]
for the optimal strategy $a^* \in\cA^{\Re}_{t,i}$, $(t,i)\in
[0,T]\times
\Ic$.
%, recalling that $C^p_L$ does not depend on $\Re$.
\end{Proposition}

%s2.3 ###
\subsection{Stability of obliquely reflected BSDEs}\label{subsecstab}
%%
%or Dependence on the parameters}

We now study the dependence on the solution with respect to the
parameters of the BSDE.
In the ``abstract'' setting considered, we obtain precious estimates
for the analysis of the regularity of the solution to the discretely
obliquely reflected BSDE as well as the convergence of the
discrete-time scheme.

We consider two discretely reflected BSDEs,
with the same reflection grid~$\Re$ but different parameters. For
$\ell\in\{ 1,2 \}$, we consider an $\Fc_T$-measurable random
terminal condition ${}^{\ell}\xi$, a random $L$-Lipschitz continuous
map\vspace*{1pt} $(y,z) \mapsto{}^{\ell}F(\cdot,y,z)$, satisfying (HF) and
a matrix of continuous cost processes\break $({}^{\ell}C^{ij})_{1\le i, j
\le d}$
satisfying the structural condition~(\ref{CostStructureabs}).

We suppose that the coefficients satisfy the integrability condition
(C$_4$).
%the following integrability condition
% \begin{equation}\label{hypoell}
% \|{}^{\ell}\xi\| _{{\bf L^2}} + \|{}^{\ell}F(\cdot,0,0)\|_{\HP{2}}
%+ \max_{i,j\in\Ic} \|{}^{\ell}C^{ij}\|_{\SP{2}} \le C_L ,
% \end{equation}
For $\ell\in\{ 1,2 \}$, we denote by
$({}^{\ell}Y^{\Re},{}^{\ell}\tY^{\Re},{}^{\ell}Z^{\Re})\in
(\mathcal{S}^{2}
\times\mathcal{S}^{2}
\times\mathcal{H}^{2})^\Ic$ the solution of the obliquely
discretely reflected BSDE
$\cD(\Re,{}^{\ell}C,{}^{\ell}F,{}^{\ell}\xi)$.\vspace*{1pt}

Defining $\dY^{\Re} = {}^1 Y^{\Re}- {}^2 Y^{\Re}$, $\dtY^{\Re} = {}^1
\tY^{\Re}- {}^2 \tY^{\Re}$, $\dZ^{\Re} =
{}^1 Z^{\Re}- {}^2 Z^{\Re}$, $\dxi:= {}^1 \xi- {}^2 \xi$ together with
\begin{eqnarray*}
|\dC_{s}|_\infty&:=& \max_{i,j\in\Ic} |{}^{1}C^{ij}-{}^{2}C^{ij}|(s)
,\\
|\dF_{s}|_\infty&:=& \max_{i \in\Ic} \sup_{y,z\in\R^d\times\Mc^{d,q}}
|{}^{1}F^{i} - {}^{2}F^{i}|(s,y,z)
\end{eqnarray*}
for $s\in[0,T]$, we prove the following stability result.
\begin{Proposition}\label{prstabY}
Assume that \textup{(C$_4$)} holds. Then we have, for any
$t\in[0,T]$,
\begin{eqnarray*}
&&
\mathbb{E}[|\dY^{\Re}_t|^2] + \mathbb{E}[|\dtY^{\Re}_t|^2]
+ \frac
{1}{\kappa} \mathbb{E}\biggl[\int_{t}^{T} |\dZ^{\Re}_{s}|^2 \,\ud s \biggr] \\
&&\qquad\le C_L
\biggl(\mathbb{E}\biggl[\int_t^T |\dF_s|^2_\infty\,\ud s + |\dxi|^2\biggr] +
\mathbb{E}\Bigl[\sup_{r\in\Re}|\dC_{r}|^4_\infty\Bigr]^{1/2}\biggr).
\end{eqnarray*}
%
%where
%|\dC_{s}|_\infty&:= \max_{i,j} |{}^{1}C^{ij}-{}^{2}C^{i,j}|(s) ,
%
%|\dF_{s}|_\infty:= \max_{i \le d} \esssup_{y,z} |{}^{1}F^{i} -
%{}^{2}F^{i}|(s,y,z).
\end{Proposition}
\begin{pf}
The proof is divided into three steps and relies heavily on the
reinterpretation in terms of switching problems. We first introduce
a convenient dominating process and then provide successively the
controls on the~$\dY^{\Re}$ and $\dZ^{\Re}$ terms.\vspace*{8pt}

\textit{Step} 1. \textit{Introduction of an auxiliary BSDE.}\vspace*{1pt}

Let us define $F:={}^{1}F\vee{}^{2}F$, $\xi:={}^{1}\xi
\vee{}^{2}\xi$ and $C$ by $C^{ij}:={}^{1}C^{ij}\wedge{}^{2}C^{ij}$.
Observe that $F$ satisfies (HF), $C$ satisfies the structure
condition~(\ref{CostStructureabs})
and that~(C$_4$) holds for the data $(C,F,\xi)$.
We denote by $(Y^{\Re},\tY^{\Re},Z^{\Re})$ the solution of the
discretely obliquely reflected
BSDE $\cD(\Re,C,F,\xi)$, recalling~(\ref{eqdeBSDEDOR1abs}).

%W_s + K _T - K _t ,\\
%K _{t}& = & \sum_{r \in\Re\setminus\set{0}} \Delta K _{r} \Ind{r
%Y_{t} &=& \tY_{t} \Ind{t \notin\Re} + \cP(t, \tY_{t})\Ind{t \in\Re}
% , t \in[0,T] .
%

Using (HF), the definition of $F$ and the monotonicity property of
$\cP$ [see Remark~\ref{repropP}(i)], we easily obtain by a comparison
argument on each interval $[r_{k},r_{k+1})$, $k\in\{0,\ldots,\kappa
-1\}
$, that
%
%e2.13 ###
%
\begin{equation}\label{eqdisccomp1}
\tY^{\Re} \succeq{}^{1}\tY^{\Re} \vee{}^{2}\tY
^{\Re}.
\end{equation}

%Since \reff{hypoell} holds,
Recalling Theorem~\ref{THswR}, we introduce the switched BSDEs
associated to ${}^{1}Y^{\Re}$, ${}^{2}Y^{\Re}$ and $Y^{\Re}$ and
denote by
$\am=(\check\theta_{j},\check a_{j})_{j\geq0}$ the optimal strategy
related to $Y^{\Re}$ starting from a fixed $(i,t)\in\Ic\times[0,T]$.
Therefore, we have
%
%e2.14 ###
%
\begin{equation}\label{eqdom0}
(\tY^{\Re}_t)^i = U^{\am}_t = \xi^{\am_T} + \int_t^T F^{\am
_s}(s,U^{\am
}_s,V^{\am}_s) \,\ud s - \int_t^T V^{\am}_s \,\ud W_s - A^{\am}_T+
A^{\am}_t.\hspace*{-28pt}
\end{equation}
%
% Following the exact same reasoning as the arguments in the proof of
%Proposition~\ref{corhatA}, we deduce from \reff{control veeY} and the
%structural cost condition that
% \begin{equation}\label{ControlNam}
% \esp{|N^{\am}|^2} \le C_L .
% \end{equation}

\textit{Step} 2. \textit{Stability of the $Y$ component.}

Since ${\am}\in\Ac^{\Re}_{t,i}$, we deduce from Theorem \ref
{THswR}(iii) that
\begin{eqnarray}
\hspace*{-3pt}&&({}^{\ell}\tY^{\Re}_t)^i \ge{}^{\ell}U^{\am_s}_{t} =
{}^{\ell
}\xi^{\am_T} + \int_t^T {}^{\ell}F^{\am_s}(s,{}^{\ell}U^{\am
}_s,{}^{\ell
}V^{\am}_s) \,\ud s - \int_t^T {}^{\ell}V^{\am}_s \,\ud W_s - {}^{\ell
}A^{\am}_T+ {}^{\ell}A^{\am}_t ,\nonumber\\
\hspace*{-3pt}&&\eqntext{\ell\in\{ 1,2 \},}
\end{eqnarray}
where ${}^{\ell}A^{\am}$ is the process of cumulated costs
$({}^{\ell}C^{ij})_{i,j\in\Ic}$ associated to the strategy~$\am$.
Combining this estimate with~(\ref{eqdisccomp1}) and~(\ref{eqdom0}), we derive
%
%e2.15 ###
%
\begin{equation}\label{eqtropfacile}
|({}^{1}\tY^{\Re}_t)^i - ({}^{2}\tY^{\Re}_t)^i | \le|U^{\am}_t -
{}^{1}U^{\am}_t| + |U^{\am}_t - {}^{2}U^{\am}_t|.
\end{equation}
Since both terms on the right-hand side of~(\ref{eqtropfacile})
are treated similarly, we focus on the first one and introduce the
continuous processes $\Gamma^{\am} := U^{\am} + A^{\am}$ and
${}^{1}\Gamma^{\am} := {}^{1}U^{\am} + {}^{1}A^{\am}$.
Applying It\^o's formula, we compute, for all $t\le u \le T$,
%
%e2.16 ###
%
\begin{eqnarray} \label{eqitostab}\quad
&&\mathbb{E}_{t} \biggl[|\Gamma^{\am}_u - {}^{1}\Gamma^{\am
}_u|^2 + \int_u^T|V^{\am}_s - {}^{1}V^{\am}_s|^2\,\ud s\biggr] \nonumber\\
&&\qquad \le\mathbb{E}_{t} \biggl[|\Gamma^{\am}_T -
{}^{1}\Gamma^{\am}_T|^2\\
&&\hspace*{48pt}{} + 2 \int_u^T (\Gamma^{\am}_s -
{}^{1}\Gamma^{\am}_s)[F^{\am_s}(s,U^{\am}_s,{}^{1}V^{\am
}_s)-{}^{1}F^{\am_s}(s,{}^{1}U^{\am}_s,{}^{1}V^{\am}_s)] \,\ud s\biggr].
\nonumber
\end{eqnarray}
Since $F={}^{1}F\vee{}^{2}F$ and ${}^{1}F$ is Lipschitz continuous,
we also get
\begin{eqnarray}
&&|F^{\am_s}(s,U^{\am}_s,{}^{1}V^{\am}_s)-{}^{1}F^{\am
_s}(s,{}^{1}U^{\am
}_s,{}^{1}V^{\am}_s)|\nonumber\\
&&\qquad \le
|\dF_{s}|_{\infty} + L ( |\Gamma^{\am}_s -{}^{1}\Gamma^{\am}_s| +
|A^{\am}_s- {}^{1}A^{\am}_s| +|V^{\am}_s - {}^{1}V^{\am}_s|) ,\qquad \eqntext{0\le s\le T.}
\end{eqnarray}
Using classical arguments, %including Gronwall Lemma,
we then deduce from the last inequality and~(\ref{eqitostab}) that
%
%e2.17 ###
%
\begin{eqnarray}\label{eqcontrolGamma}
|\Gamma^{\am}_t\! -\! {}^{1}\Gamma^{\am}_t|^2\!\le\! C_L \biggl(
\mathbb{E}_{t} \biggl[|\dxi^{\am_T}|^2 \int_t^T |\dF
_{s}|_{\infty}^{2} \,\ud s\biggr]\! +\! \sup
_{t\le s
\!\le\! T} \mathbb{E}_{t} [|A^{\am}_s\!-\! {}^{1}A^{\am}_s|^2
] \biggr).\hspace*{-40pt}
\end{eqnarray}
Moreover, using the inequality $|x\vee y - y|\leq|x-y|$ for
$x,y\in\R$ and the convexity of the function $x\mapsto x^2$, we compute
%
%e2.18 ###
%
\begin{eqnarray} \label{eqcontrolA}
&&\mathbb{E}_{t} [|A^{\am}_s- {}^{1}A^{\am}_s|^2] \nonumber\\
&&\qquad=
\mathbb{E}_{t} \Biggl[\Biggl| \sum_{k= 1}^{\Nm}
[^2C^{\alfam_{k-1}\alfam_{k}}\wedge^1C^{\alfam_{k-1}\alfam
_{k}}-^1C^{\alfam_{k-1}\alfam_{k}} ](\tetam_{k}) \mathbf
{1}_{ \{ \tetam_{k} \le s \} } \Biggr|^2 \Biggr]\\
&&\qquad\le \mathbb{E}_{t} \Bigl[|N^{\am}| \sup_{r \in\Re}|\dC
_{r}|_{\infty}^2\Bigr],\qquad
t \le s \le T.\nonumber
\end{eqnarray}
Plugging \label{eqcontrolA} in~(\ref{eqcontrolGamma}) and recalling
the definition of $\Gamma^{\am}$ and ${}^{1}\Gamma^{\am}$, we get
\[
{|U^{\am}_t - {}^{1}U^{\am}_t|^2} \le C_L\mathbb{E}_{t}
\biggl[|N^{\am}| \sup_{r \in\Re}|\dC_{r}|_{\infty}^2 + {\int
_t^T |\dF_{s}|_{\infty}^{2} \,\ud s + |\dxi|^2}\biggr].
\]
The exact same reasoning leads to the same estimate for ${|U^{\am}_t
- {}^{2}U^{\am}_t|^2}$. Therefore, we deduce from~(\ref{eqtropfacile})
and the Cauchy--Schwarz inequality that
%
%e2.19 ###
%
\begin{eqnarray}\label{eqcontrolstabYCOR}\qquad
&&\mathbb{E}[|({}^{2}\tY^{\Re}_t)^i -
({}^{1}\tY^{\Re}_t)^i|^2]\nonumber\\[-8pt]\\[-8pt]
&&\qquad\le C_L
\biggl(\mathbb{E}[|N^{\am}|^2]^{1/2} \mathbb{E}\Bigl[\sup_{r \in\Re
}|\dC_{r}|_{\infty}^4\Bigr]^{1/2} + \mathbb{E}\biggl[{\int_t^T |\dF
_{s}|_{\infty}^{2} \,\ud s + |\dxi|^2}\biggr]\biggr).\hspace*{-20pt}\nonumber
\end{eqnarray}
%
% \esp{|\dtY_t|^2} \le C_L ( \esp{ \int_t^T |\dF_s|^2_\infty\ud
%s + |\dxi|^2} +\esp{\sup_{r\in\Re}|\dC_{r}|^4_\infty}^{1 \over2} )
% .
Using Proposition~\ref{corhatA}, we compute, since $i$ is
arbitrary,
%
%e2.20 ###
%
\begin{eqnarray}\label{eqcontrolstabY}
&&
\mathbb{E}[|{}^{2}\tY^{\Re}_t - {}^{1} \tY^{\Re}_t|^2] \nonumber\\[-8pt]\\[-8pt]
&&\qquad\le C_L
\biggl(
\mathbb{E}\biggl[\int_t^T |\dF_s|^2_\infty\,\ud s + |\dxi|^2\biggr] +\mathbb
{E}\Bigl[\sup_{r\in\Re}|\dC_{r}|^4_\infty\Bigr]^{1/2} \biggr).\nonumber
\end{eqnarray}

%
%M_{t} := \EFp{t}{ \beta\sup_{r \in\Re}|\dC_{r}|_{\infty}^2 + {
%we have that
%Using then Burkholder-Davis-Gundy inequality, we compute
%which concludes the proof of the proposition.
%Using Proposition~\ref{prstabY}, we are able to obtain a control on
%the variation of the $Z$-component of the discretely obliquely
%reflected BSDEs
%The following holds,
%+ |\dxi|^2 } + \esp{\sup_{r \in\Re}|\dC_{r}|_{\infty}^4}^\frac12 )

\textit{Step} 3. \textit{Stability of the $Z$ component.}

Applying It\^o's formula to the \textit{c\`{a}dl\`{a}g} process $|\dtY
^{\Re}|^{2}$
and noting $\dtK= {}^{1}K ^{\Re} -{}^{2}K ^{\Re} $, we obtain
\begin{eqnarray*}
&&\mathbb{E}\biggl[|\dtY^{\Re}_{t} |^{2} + \int_{t}^{T} |\dZ^{\Re
}_{s}|^2 \,\ud s + \sum_{t< r \le T}|\Delta\dtK^{\Re}_{r}|^2\biggr]\\
&&\qquad =
\mathbb{E}\biggl[|\dtY^{\Re}_{T}|^2 + 2 \int_{t}^T \dY^{\Re}_{s}
\dF_{s} \,\ud s + 2 \int_{t}^T \dY^{\Re}_{r} \,\ud\dtK^{\Re
}_{r} \biggr],
\end{eqnarray*}
where we used the fact that $|\dtY^{\Re}|^2 - |\dY^{\Re}|^2 - 2\dY
^{\Re
}(\dtY^{\Re}- \dY^{\Re}) = |\Delta\dtK^{\Re}|^2$.
Since $\delta K$ is a pure jump process, we compute
\[
\mathbb{E}\biggl[\int_{t}^T \dY^{\Re}_{r} \,\ud\dtK^{\Re}_{r} \biggr] \le
\mathbb{E}\biggl[\alpha\sum_{t <r \le T, r\in\Re} |\dY^{\Re}_{r}|^2 +
\frac{1}{\alpha} \sum_{t< r \le T}|\Delta\dtK^{\Re}_{r}|^2
\biggr],\qquad
\alpha>0,
\]
which, for $\alpha$ large enough and using standard arguments, leads to
\begin{eqnarray*}
&&\mathbb{E}\biggl[\int_{t}^{T} |\dZ^{\Re}_{s}|^2 \,\ud s + \sum_{t<
r \le T}|\Delta\dtK^{\Re}_{r}|^2 \biggr]\\
&&\qquad \le C_L \biggl( \mathbb
{E}[|\dxi|^2] +
\mathbb{E}\biggl[\int_t^T |\dF_{s}|^2_\infty\,\ud s + \sum_{t <r
\le T,r\in\Re} |\dY^{\Re}_{r}|^2\biggr] \biggr).
\end{eqnarray*}
Since~(\ref{eqcontrolstabY}) holds true for any $t\in[0,T]$, we deduce
\begin{eqnarray*}
&&\mathbb{E}\biggl[\int_{t}^{T} |\dZ^{\Re}_{s}|^2 \,\ud s + \sum_{t<
r \le T}|\Delta\dtK^{\Re}_{r}|^2 \biggr]\\
&&\qquad \le C_L \kappa\biggl( \mathbb
{E}[|\dxi|^2]
+ \mathbb{E}\biggl[\int_t^T |\dF_{s}|^2_\infty\,\ud s \biggr] + \mathbb
{E}\Bigl[\sup_{r \in\Re}|\dC_{r}|_{\infty}^4\Bigr]^{1/2} \biggr),
\end{eqnarray*}
which concludes the proof of the proposition.
\end{pf}

%%%%%%%%%%%%%%%%%%%%%%
%%%%%%%%%%%%%%%%%%%%%%
%%%%%%%%%%%%%%%%%%%%%%
%%%%%%%%%%%%%%%%%%%%%%

%!TEX root = main.tex

%s3 ###
\section{Regularity of discretely obliquely reflected BSDEs}\label
{secreg}

This section is dedicated to the derivation of regularity properties
for the solution of discretely reflected BSDEs. These results are
obtained in a Markovian diffusion setting. This means that the randomness
of the parameter $(C,F,\xi)$ is due to a state process $X$, which
is the solution of a stochastic differential equation (SDE). In this
framework, we focus on the $\Hc^2$-regularity of the $Z^{\Re}$ component
of the solution of the BSDEs. The main results are retrieved by
means of kernel regularization and Malliavin differentiation
arguments. Finally, we extend this result to the case where the
diffusion $X$ is replaced by its Euler scheme.

%s3.1 ###
\subsection{A diffusion setting for discretely RBSDEs}

Let $X$ be the solution on $[0,T]$ to the
following SDE:
%
%e3.1 ###
%
\begin{equation}\label{EqForward}
X_t = X_{0} + \int_0^t b(X_s) \,\ud s
+ \int_0^t
\sigma(X_s) \,\ud W_s ,\qquad 0\le t\le T,
\end{equation}
where $X_{0} \in
\R^m$ and $(b,\sigma)\dvtx\R^m\rightarrow\R^m\times\Mc^{m,q}(\R)$ are
$L$-Lipschitz functions.

%and introduce the process $X$, unique solution of the following
%stochastic differential equation

Under the above assumption, the following estimates are well known
(see, e.g.,~\cite{kun90}):
%
%e3.2 ###
%
\begin{eqnarray}\label{eqcontrolX}
\mathbb{E}\Bigl[\sup_{t \in
[0,T]}|X_t|^p\Bigr] &\le& C_L^p \quad\mbox{and}\nonumber\\[-8pt]\\[-8pt]
\sup_{s\in[0,T]}
\Bigl(\mathbb{E}\Bigl[\sup_{u\in[0,T],|u-s|\le h}|X_{s}-X_{u}|^p\Bigr] \Bigr)
^{1/p} &\le&
C^p_L \sqrt{h}\nonumber
\end{eqnarray}
%
%}\end{Remark}
for any $p>0$. In the sequel, we shall denote by $\beta^X$ a positive
random variable, which may change from line to line,
but which depends only on $\sup_{t \in[0,T]} |X_{t}|$ and which
satisfies $\mathbb{E}[|\beta^X|^p]\le C^p_{L}$ for all $p > 0$. Importantly,
$\beta^X$ does not depend on~$\Re$.
\begin{Remark}
Observe that, as in~\cite{boucha08,cha09} and contrary to
\cite{majzha05}, we make no uniform ellipticity condition on
$\sigma$. This allows us to treat the case of nonhomogenous
diffusion by setting, for example, $X^1_t = t$, $t \in[0,T]$.
\end{Remark}

In this context, we are given a matrix valued
maps $c:= (c^{ij})$ where
$c^{ij}\dvtx\R^m\rightarrow\R^+$, are $L$-Lipschitz continuous and satisfy
%
%e3.3 ###
%
\begin{equation}\label{CostStructure}\qquad
\cases{
\displaystyle c^{ii}(\cdot)=0 , &\quad for $1\le i \le d$; \cr
\displaystyle \inf_{x\in\R^m}c^{ij}(x)>0 , &\quad for $1\le i,j \le d$ with $i\neq j$;
\cr
\displaystyle
\inf_{x\in\R^m}\{c^{ij}(x) + c^{jl}(x) - c^{il}(x)\}>0, &\quad for
$1\le i,j,l\le d$\cr
&\quad
with $i\neq j, j\neq l$.}
\end{equation}

We then introduce a family $(\cQ(x))_{x \in\R^m}$ of closed convex
domains:
%
%e3.4 ###
%
\begin{eqnarray}\label{eqdefcC}
\cQ(x) := \Bigl\{ y \in\R^d
\bigm| y^i
\geq\max_{j\in\Ic} \bigl(y^j - c^{ij}(x)\bigr) , \forall i \in\cI
\Bigr\}\nonumber\\[-8pt]\\[-8pt]
&&\eqntext{\mbox{where } \cI:=\{ 1,\ldots,d \}.}
\end{eqnarray}

%conditions allow to get existence and uniqueness of a solution to the
%corresponding continuously reflected BSDE, see Section
%Moreover, they are economically meaningful
We introduce the oblique projection operator $\cP(x,\cdot)$ onto
$\cQ(x)$ defined by
\[
\cP\dvtx(x,y) \in\R^m\times\R^d \mapsto\Bigl( \max_{j\in\Ic
}
\{ y^j - c^{ij}(x) \} \Bigr)_{1\le i\le d}.
\]

%It follows from the structure condition \reff{Cost Structure} that $
% \begin{equation}\label{PLipschitz}
% | \cP(x_{1},y_1) - \cP(x_{2},y_2) |
% &= ( \sum_{i=1}^d | \Max_{j\le d} (y^j_1-c_{ij}(x_{1})) -
% &\le\sqrt{d} ( \Max_{j\le d} |y_1^j - y_2^j|+ \Max_{i,j\le
%d}|c_{ij}(x_1) - c_{ij}(x_2)| ) \nonumber\\
% & \le L (|y_1 - y_2 |+|x_1 - x_2 |) ,
% \end{equation}
% for all $ y_1, y_2 \in\R^d$, and $x_1, x_2 \in\R^m$. Due to the
%obliqueness of the projections, this operator is not $1$-Lipschitz,
%ruling out classical methods for the analysis conducted in this paper.

Finally, we are given:

\begin{longlist}
\item
an $L$-Lipschitz function $g\dvtx\R^m\rightarrow\R^d$ such
that $g(x)\in\cQ(x)$ for all $x\in\R^m$,%a terminal condition, i.e.
%an $

\item a generator function, that is, an $L$-Lipschitz map $f\dvtx
\R^m\times\R^d\times\Mc^{d,q}\rightarrow\R^d$.
\end{longlist}

From now on, we shall appeal to the following assumption:

\begin{longlist}[(Hf)]
\item[(Hf)] the component $i$ of $f(\cdot,y,z)$ depends only on the
component $i$ of the vector $y$ and on the column $i$ of the matrix
$z$, that is, $f^i(\cdot,y,z)=f^i(\cdot,y^i,z^i)$.
\end{longlist}

We denote by $(Y^{\Re},\tY^{\Re},Z^{\Re})$ the solution of the
discretely reflected BSDE $\cD(\Re, c(X), f(X,\cdot,\cdot), g(X))$ which reads
on each interval $[r_j,r_{j+1})$, for $j<\kappa$
%
%e3.5 ###
%
\begin{equation}\label{eqdeBSDEDOR1}
\cases{
\displaystyle \tY^{\Re}_{t } = Y^{\Re}_\rjp+ \int_{t}^\rjp
f(X_{u}, \tY^{\Re}_u, Z^{\Re}_u)\,\ud u - \int_{t}^\rjp Z^{\Re}_u
\,\ud W_u,
\vspace*{2pt}\cr
\displaystyle Y^{\Re}_{t} = \tY^{\Re}_t\mathbf{1}_{ \{ t \notin\Re\} } + \cP
(X_{t},\tY^{\Re
}_t)\mathbf{1}_{ \{ t \in\Re\} } ,}
\end{equation}
or equivalently on $[0,T]$ as
%
%e3.6 ###
%
\begin{equation}\label{eqBSDEDOR}
\cases{
\displaystyle \tY^{\Re}_{t } = g(X_{T}) + \int_{t}^T f(X_{u},\tY^{\Re
}_u,Z^{\Re}_u)\,\ud u \vspace*{2pt}\cr
\hphantom{\displaystyle \tY^{\Re}_{t } =}{}- \displaystyle \int_{t}^T Z^{\Re}_u \,\ud W_u + (K ^{\Re}_T-K
^{\Re}_t), &\quad$0\le t\le T$, \vspace*{2pt}\cr
\displaystyle K ^{\Re}_t :=
\sum_{r \in\Re\setminus\{ 0 \}} \Delta K ^{\Re}_r \mathbf{1}_{
\{ r \le t \} }\quad
\mbox{and }\vspace*{2pt}\cr
\displaystyle \Delta K ^{\Re}_t = Y ^{\Re}_t - \tY^{\Re}_t = -(\tY^{\Re}_t -
\tY^{\Re
}_{t-}), &\quad$0\le t \le T$.}
\end{equation}

%Observe that $Y$ and $\tY$ differ only on the grid points of $\Re$.
%On each interval of the form $[\rk,\rkp)$, $(\tY,Z)$ is solution to a
%classical non reflected BSDE with terminal condition $Y_{\rkp}$, see
%the discretely reflected BSDE \eqref{eqdeBSDEDOR1} follows
%directly from a concatenation of the solutions on all the grid
%intervals. We finally provide some a priori estimates on the solution
%to \eqref{eqdeBSDEDOR1}.% , whose proof is reported in the
%Appendix~\ref{AppendixAprioriDOR} of the paper.

From~(\ref{eqcontrolX}), it follows that the data $(c(X),
f(X,\cdot,\cdot), g(X))$ satisfies the integrability condition~(\ref{Cp})
for all $p \ge2$. We thus deduce from the proofs of Propositions
\ref{estprioedsrdrabs} and~\ref{corhatA}, the
following estimate on $(Y^{\Re},\tY^{\Re},Z^{\Re})$ and their
associated optimal
switched BSDEs, recalling Theorem~\ref{THswR}.
\begin{Proposition}\label{estprioedsrdr}
There exists a unique solution $(\tY^{\Re},Y^{\Re},Z^{\Re})$ to
(\ref{eqdeBSDEDOR1}) and it satisfies
%
%e3.7 ###
%
\begin{equation}\label{controlContRef}\qquad
\mathbb{E}_{t} \biggl[\sup_{ s\in[t,T]} |\tY^{\Re}_{s}|^p +
\biggl(\int_{t}^T|Z_{s}^{\Re}|^2\,\ud s\biggr)^{p/2} + |K^{\Re
}_{T}-K^{\Re}_{t}|^p \biggr] \le\beta
^X\qquad \forall t \le T.\hspace*{-28pt}
\end{equation}
Moreover, for all $(t,i)\in[0,T]\times\cI$, the optimal strategy $a^*
\in\cA^\Re_{t,i}$ satisfies
%
%e3.8 ###
%
\begin{equation}\qquad
\mathbb{E}_{t} \biggl[\sup_{s\in[t,T]}|U^{a^*}_{s}|^p +\biggl(\int
_{t}^T|V^{a^*}_{s}|^2\,\ud s\biggr)^{p/2}+ |A^{a^*}_{T}|^p + |N^{a^*}|^p
\biggr]
\le\beta^X.
\end{equation}
%
%for all $(t,i)\in[0,T]\times\Ic$.
\end{Proposition}

% \romu{Discuter comment relacher cette hypothese \jef{regularisation
%par noyau}}

% \jef{definir les notations pour les derivees...}

%s3.2 ###
\subsection{\texorpdfstring{Malliavin differentiability of $(X,Y^{\Re},\tY^{\Re},Z^{\Re})$}
{Malliavin differentiability of (X, Y^R, Y^R, Z^R)}}

We shall sometimes use the following regularity assumption on the coefficients:

\begin{longlist}[(Hr)]
\item[(Hr)] The coefficients $b$, $\sigma$, $g$ $f$ and $(c^{ij})_{i,j}$
are $C^{1,b}$ in all their variables, with the Lipschitz constants
dominated by $L$.
\end{longlist}

We denote\vspace*{-1pt} by $\D$ the set of random variables $G$ which are
differentiable in the Malliavin sense and such that ${\|G\|_\D}^2:=\|
G\|
^2_{\Lb^2} + \int_0^T \|D_t G\|_{\Lb^2}^2\,\ud t<\infty$, where $D_t G$
denotes the Malliavin derivative of $G$ at time $t\le T$.
After\vspace*{1pt}
possibly passing to a suitable version, an adapted process belongs to
the subspace $\Lc_a^{1,2}$ of $\Hc^2$ whenever $V_s\in\D$ for all
$s\le
T$ and $\|V\|_{\Lc_a^{1,2}}^2:=\|V\|^2_{\Hc^2} + \int_0^T \|D_t V\|
_{\Hc
^2}^2\,\ud t<\infty$. For\vspace*{-1pt} a general presentation on Malliavin
calculus for stochastic differential equations, the reader may refer to
\cite{Nua95}.
\begin{Remark}\label{RemarkMalliavinDiffX}
Under (Hr), the solution of~(\ref{EqForward}) is Malliavin
differentiable and its derivative satisfies
%
%e3.9 ###
%
\begin{equation} \label{eqmajoDX}
\Bigl\| {\sup_{s\le T} }| D_s X | \Bigr\|_{_{\Sc^p}} < \infty,
\end{equation}
and we have
%
%e3.10 ###
%
\begin{eqnarray}\label{ControlReguDsXt}
&&\sup_{s\le u } \| D_sX_t - D_s X_u \|_{_{\mathbf{L}^p}} +
\Bigl\|{\sup_{t \le s \le T }} | D_tX_s - D_u X_s |
\Bigr\|_{_{\mathbf{L}^p}}\nonumber\\[-9pt]\\[-9pt]
&&\qquad\le C_L^p | t-u |^{1/2}\nonumber
\end{eqnarray}
for any\vspace*{1pt} $0\le u\le t\le T$. Let $G \in\D(\R^d)$. Since $X$ belongs to
$\Lc_a^{1,2}$ under (Hr) and $\cP$ is $L_\cP$-Lipschitz continuous,
%according to \eqref{PLipschitz},
we deduce that $\cP(X_{t},G) \in\D(\R^{d})$. Using Lemma 5.1 in
\cite
{boucha08}, we compute
%
%e3.11 ###
%
\begin{eqnarray}
\label{eqmalderivcP}
&& D_{s}(\cP(X_{t},G))^{i} \nonumber\\[-2pt]
&&\qquad= \sum_{j = 1}^{d} \bigl(D_{s}G^{j} - D_{s}c_{ij}(X_{t})\bigr)\mathbf{1}_{ \{
G^{j} - c^{ij}(X_{t}) > \max_{\ell<j}(G^{\ell}-c^{i \ell}(X_{t}))
\} } \\[-2pt]
&&\qquad\quad\hspace*{13pt}{}\times\mathbf{1}_{ \{ G^{j} - c^{ij}(X_{t}) \geq\max_{\ell
>j}(G^{\ell}-c^{i \ell}(X_{t})) \} } .
\nonumber
\end{eqnarray}
\end{Remark}

Combining~(\ref{eqmalderivcP}), Proposition 5.3 in \cite
{elkkapparpenque97} and an induction argument, we obtain that
$(Y^{\Re},\tY^{\Re},Z^{\Re})$ is Malliavin differentiable and that a
version of $(D_{u} \tY^{\Re}, D_{u}Z^{\Re})$ is given by
% \begin{Lemma}\label{LemmaDY}
% Under \HYP{r}, the solution $(Y,\tY,Z)$ of \reff{eq de BSDE DOR 1} is
%in $\Lc_a^{1,2}$ and its Malliavin %derivative $(D_uY,D_u\tY,D_uZ)$,
%for $u\le T$, is solution on each grid interval $[r_j,r_{j+1})$ of the
%%following linear BSDE
%
%e3.12 ###
%
\begin{eqnarray}\label{EqDY}
D_u (\tY^{\Re}_t)^i &=& D_u (Y^{\Re}_{\rjp})^i
- \sum_{k=1}^d \int_t^{r_{j+1}} D_u (Z^{\Re}_s)^{ik} \,\ud W^k_s\nonumber\\[-2pt]
&&{} + \int_t^{\rjp} \nabla_{x} f^i (X_s,(\tY^{\Re
}_s)^i,(Z^{\Re}_s)^{i\cdot}) D_uX_s\,\ud s \nonumber\\[-9pt]\\[-9pt]
&&{}+ \int_t^{\rjp} \nabla_{y^i} f^i (X_s,(\tY^{\Re
}_s)^i,(Z^{\Re}_s)^{i\cdot}) D_u(\tY^{\Re}_s)^i \,\ud s\nonumber\\[-2pt]
&&{}+ \int_t^{\rjp} \nabla_{z} f^i (X_s,(\tY^{\Re}_s)^i,(Z^{\Re}_s)^{i\cdot})
D_u(Z^{\Re}_s)^{i\cdot} \,\ud s\nonumber
% - \sum_{k=1}^d \int_t^{r_{j+1}} D_u (Z^{\Re}_s)^{ik} \,\ud W^k_s ,
\end{eqnarray}
for $0\le u\le t\le r_{j+1}$ and $j<\kappa$. Here, $\nabla_{z}f^{i}$
denotes $\sum_{\ell= 1}^{d}\nabla_{z^{\ell.}} f^i $, recalling~(Hf).

%{\bf Proof.} The Malliavin differentiability of $(\tY,Y,Z)$ follows
%from an induction argument.
%First, the terminal value $Y_T=\tY_T$ is in $\D$ since $g$ is
%$C^{1,b}$ and $X\in\D$.
% Recursively, fix $k/in{\kappa,\ldots,1}$ and suppose $\tY_{r_k}\in\D$.
% Observe that $Y_{r_k}\in\D$ since the projection operator $\Pc$ is
%Lipschitz, see Lemma~\ref{LemmaP}.
%Following similar arguments as the one of Proposition in 5.3. in

%The behavior of the Malliavin derivatives of $(\tY,Y)$ on the grid
%points is detailed in the next Lemma.

% \begin{Lemma}\label{LemmaDYgrid}
% Fix $(\theta_0,\alpha_0)\in[0,T]\times\Ic$ and introduce the optimal
%switching strategy $a$ associated to $\tY$ as in Theorem~\ref{THswR}.
%Under \Ha, we have
% \be
% D_t Y^{\alpha_{j-1}}_{\theta_j} = D_t \tY^{\alpha_{j}}_{\theta_j} -
%c_{\alpha_{j-1},\alpha_j}(X_{\theta_j}) \{D_tX\}_{\theta_j} .
% \ee
%D_t Y^{a_{r_j}}_{r_j} = D_t \tY^{a_{r_j+}}_{r_{j}} -
%c'_{a_{r_j},a_{r_j+}}(X_{r_j}) D_tX_{r_j} , 0\le j<\kappa,
% t\le T .

%{\bf Proof.} At a grid point where $a$ does not jumps, we have $Y^a=
%At a jump time $\theta_j$ where $a$ jumps from $\alpha_{j-1}$ to $

%s3.3 ###
\subsection{Representation of $Z$}% and regularity of $Z$}
\label{sserepz}

For $a \in\cA^{\Re}$, we introduce the process $\Lambda^{a}$
defined by
%
%e3.13 ###
%
\begin{eqnarray}\label{defLambda}
\Lambda^{a}_{t,s} &:=& \exp\biggl\{ \int_t^s \nabla_{z}
f^{a_r}(X_r,\tY
^{\Re}_r,Z^{\Re}_r) \,\ud W_r \nonumber\\[-9pt]\\[-9pt]
&&\hphantom{\exp\biggl\{}
{} - \int_t^s \biggl(
\frac
{1}{2} |\nabla_{z} f^{a_r}(X_r,\tY^{\Re}_r,Z^{\Re}_r)|^2 -
\nabla
_{y}f^{a_r}(X_r,\tY^{\Re}_r,Z^{\Re}_r) \biggr) \ud r \biggr\}\hspace*{-28pt}\nonumber
\end{eqnarray}
for $0\le t\le s\le T$.

For later use, we remark
%
%e3.14 ###
%
\begin{equation}\label{eqlambdasup}
\sup_{a \in\cA^{\Re}} \Bigl\| \sup_{t\le s\le T} \Lambda_{t,s}^{a}\Bigr\|
_{\cL^p
}\le C^p_L ,\qquad 0\leq t \leq T, p\ge2,\vadjust{\goodbreak}
\end{equation}
and deduce from the dynamics of $\Lambda$ that
%
%e3.15 ###
%
\begin{eqnarray}\label{eqlambdacontrol1}
&&
\sup_{a \in\cA^{\Re}} \Bigl( \| \Lambda^{a}_{t,t} - \Lambda^{a}_{t,u}
\|_{\Lc^p} + \Bigl\|\sup_{t \le s \le T } | \Lambda^a_{u,s} - \Lambda
^a_{t,s} | \Bigr\|_{\Lc^p} \Bigr)\nonumber\\[-8pt]\\[-8pt]
&&\qquad\le C^p_L \sqrt{t-u} ,\qquad u\le t\le T , p\ge2 .\nonumber
%& \le C^p_L \sqrt{t-s} , {a \in\cA^{\Re}_{s,i}} ,{b \in\cA^{
\end{eqnarray}

% which combined with \reff{Control Regu DsXt} leads to
% \be\label{ControlReguDsXt2}
% \sup_{a_t\in\Ic} \| \sup_{s \le u \le T} | \Lambda_u^{t,a_t} D_tX_u
%- \Lambda^{s,a_s}_u D_s X_u | \|_{\Lc^p} & \le& C_L^p |t-s|^{
% \ee
%, where $a$ is the optimal strategy associated to $(Y,\tY,Z)$
%starting from $a_t$ at time $t$. A particular %feature of Malliavin
%differentiable BSDE is the reinterpretation of the process $Z$ as a
%version of
%$(D_tY_t)_{t\le T}$. This remarkable property combined with the
%previous lemmata lead to the following %representation for $Z$.
%
\begin{Proposition}\label{PropRepresZ}
Under \textup{(Hr)}, there is a version of $Z^{\Re}$ such that
%
%e3.16 ###
%
\begin{eqnarray}\label{EqRepresZ}
(Z^{\Re}_t)^i &=& \mathbb{E}_{t} \Biggl[\nabla_{x}g^{a^{*}_T}(X_T)
\Lambda^{a^{*}}_{t,T} D_tX_T\nonumber\\
&&\hphantom{\mathbb{E}_{t} \Biggl[}
{} + \int_t^T \nabla
_{x} f^{a^{*}_s}(X_s,\tY^{\Re}_s,Z^{\Re}_s) \Lambda^{a^{*}}_{t,s}
D_tX_s \,\ud s \\
&&\hphantom{\mathbb{E}_{t} \Biggl[}\hspace*{15pt}
{} -\sum_{j=1}^{N^{a^{*}}}
\nabla_{x} c^{\alpha^{*}_{j-1}\alpha^{*}_j}(X_{\theta^{*}_j})
\Lambda^{a^{*}}_{t,{\theta^{*}_j}}(D_t X)_{\theta^{*}_j}\Biggr]\nonumber
\end{eqnarray}
for $(t,i)\in[0,T]$, with $a^{*}=(\theta^*_j,\alpha^*_j)_{j\ge0}
\in
\cA^{\Re}_{t,i}$ the optimal strategy given in Theorem~\ref{THswR}
and recalling~(\ref{eqdefNa}).
\end{Proposition}
%
%theta j-1, theta j}
%
\begin{pf}
We fix $j < \kappa$ and, observing that the process $a^*$ is
constant on the interval $[\theta_{j}^*,\theta_{j+1}^*)$, we deduce
from~(\ref{EqDY}) and It\^o's formula that
\begin{eqnarray*}
&&\Lambda^{a^*}_{t,t} D_u(\tY^{\Re}_t)^{\alpha^*_j} \\
&&\qquad=
\E_{t}\biggl[ \Lambda^{a^*}_{t,\theta^*_{j+1}} (D_u (Y^{\Re
})^{\alpha
^*_j})_{\theta^*_{j+1}} + \int_t^{\theta^*_{j+1}} \nabla_x
f^{\alpha
^*_j} (X_s,\tY^{\Re}_s,Z^{\Re}_s) \Lambda^{a^*}_{t,s} D_uX_s \,\ud s
\biggr]
% \\ & -\int_t^{T} (D_uZ)^{a_s}_s\cdot dW_s
% - \sum_{j=1}^{N^{a^{*}}} \nabla c_{ \alpha^{*}_{j-1} , \alpha^{*}_j }
%(X_{\theta^{*}_j}) (D_uX)_{ \theta^{*}_j} ,
\end{eqnarray*}
for\vspace*{1pt} $\theta^*_j \le u \le t < \theta^*_{j+1}$.
Combining~(\ref{eqmalderivcP}) and the definition of $a^{*}$ given
in Theorem~\ref{THswR}(ii), we compute
\begin{eqnarray*}
\Lambda^{a^*}_{t,\theta^*_{j+1}} (D_u (Y^{\Re})^{\alpha
^*_j})_{\theta^*_{j+1}}
&=& \Lambda^{a^*}_{t,\theta^*_{j+1}} (D_u (\tY^{\Re})^{\alpha
^*_{j+1}})_{\theta^*_{j+1}}\\
&&{}- \nabla_{x} c^{\alpha^{*}_{j}\alpha^{*}_{j+1}} (X_{\theta
^{*}_{j+1}}) \Lambda^{a^{*}}_{t,\theta^*_{j+1}} (D_t X)_{\theta^{*}_{j+1}}
\end{eqnarray*}
for $j < \kappa$. Plugging the second equality into the first one and
summing up over $j$ concludes the proof.
\end{pf}
%

%definition of $a^{*}$ given in Theorem~\ref{THswR} (ii) and iterating
%over the points of $\Re$ greater than $t$, we get
% \[
% D_u(\tY_t)^{i} &= (D_u Y)^{a_T}_T + \int_t^{T} (\nabla f^{a_s} (X_s,
% \\ & -\int_t^{T} (D_uZ)^{a_s}_s\cdot dW_s
% - \sum_{j=1}^{N^{a^{*}}} \nabla c_{ \alpha^{*}_{j-1} , \alpha^{*}_j }
%(X_{\theta^{*}_j}) (D_uX)_{ \theta^{*}_j} ,
% \]
% recall that $c_{jj}=0$, $j \in\cI$.
% Observe that the last term of the previous expression coincides with
%the last term of \reff{Eq Repres Z} since $c_{j,j}=0$, for $j\le
% Applying Ito's formula to the process $\Lambda^{a^{*}}D_u(\tY)^{i}$
%between $t$ and $T$, conditioning with respect to $\Fc_t$, %and
%observing that $(c'_{jj})_{j\le\kappa}=0$ and $D_uX_{\theta_j}=0$ for
%$u\ge\theta_j$
%we derive \reff{Eq Repres Z}.

We conclude this section by providing a ``weak'' regularity property of~$Z^{\Re}$ in the general Lipschitz setting.
In order to get rid of the previous assumption~(Hr), we make use of
kernel regularization arguments.
Since this procedure is very classical, we do not detail it here
precisely (see, e.g., the proofs of Proposition 4.2 in~\cite{cha09} or
Proposition 3.3 in~\cite{boucha08}).
%We omit to detail precisely this classical procedure, which is not the
%main concern of the paper.
%
\begin{Proposition}\label{prZbounded}
There is a version of $Z^{\Re}$ satisfying
%
%e3.17 ###
%
\begin{equation}\label{ReguweakZ}
\mathbb{E}\biggl[\int_{s}^t |Z^{\Re}_{u}|^{2} \,\ud u\biggr]
\le C_{L} |t-s| ,\qquad
s\le
t \le T.\vadjust{\goodbreak}
\end{equation}
%
%(ii) Moreover, under \HYP{r},
% \[
\end{Proposition}
\begin{pf}
Combining~(\ref{eqmajoDX}), with~(\ref{eqlambdasup}),~(\ref{EqRepresZ})
and Doob's inequality, we observe that
\[
\sup_{t\in[0,T]} \|Z^{\Re}_{t}\|_{_{\mathbf{L}^p}}\leq C_L^p,\qquad
p\ge2,
\]
holds under (Hr). Therefore~(\ref{ReguweakZ}) is satisfied under
(Hr).
As in the proof of Proposition 4.2 in~\cite{cha09}, the stability
results of Proposition~\ref{prstabY} allow us to use classical Kernel
regularization arguments. Since the previous estimate holds uniformly
for the sequence of approximating regularized BSDE, the proof is complete.
\end{pf}

%boucha08 aussi}

%s3.4 ###
\subsection{Regularity results}
\label{sseregz}

We consider a grid $\pi:=\{t_0=0,\ldots,t_n=T\}$ on the time
interval $[0,T]$, with modulus $|\pi|:=\max_{0\leq i\leq
n-1}|t_{i+1}-t_{i}|$, such that $\Re\subset\pi$.
% and $|\pi|n\leq L$.

We want to control the following quantities, representing the $\mathcal
{H}^{2}$-regularity of $(\tY,Z)$:
%
%e3.18 ###
%
\begin{equation}\quad
\mathbb{E}\biggl[\int_0^T\bigl|\tY^{\Re}_t - \tY^{\Re}_{\pi( t)}\bigr|
^2 \,\ud t\biggr]
\quad\mbox{and}\quad \mathbb{E}\biggl[\int_0^T
\bigl| Z^{\Re}_t - \bar{Z}^{\Re}_{\pi
(t)} \bigr|^2 \,\ud t\biggr],
\end{equation}
where $\pi(t):=\sup\{\ti\in\pi; \ti\le t\}$ is defined on $[0,T]$
as the projection to the closest previous grid point of $\pi$
and
%
%e3.19 ###
%
\begin{equation}\label{defbarZ}\quad
\bar Z^{\Re}_\ti:= \frac{1}{\tip-\ti} \mathbb{E}\biggl[\int_\ti
^\tip Z^{\Re}_s \,\ud s\Bigm|\mathcal{F}_{\ti}\biggr] ,\qquad
i\in\{0,\ldots,n-1 \}.
\end{equation}

\begin{Remark} \label{rebestapprox}
Observe that $(\bar Z^{\Re}_s)_{ s\le T}:=(\bar Z^{\Re}_{\pi(s)})_{ s
\le T}$ is interpreted as the best $\mathcal{H}^{2}$-approximation of the
process $Z^{\Re}$ by adapted processes which are constant on each
interval $[\ti,\tip)$, for all $i<n$.
\end{Remark}
\begin{Proposition} \label{prregY} The following holds:
\[
\frac{1}{T} \mathbb{E}\biggl[\int_0^T\bigl|\tY^{\Re}_t - \tY^{\Re}_{\pi
( t)}\bigr|^2 \,\ud t\biggr] \le\sup_{t \in[0,T]}
\mathbb{E}\bigl[\bigl|\tY^{\Re}_t -
\tY^{\Re}_{\pi(t)}\bigr|^2 \bigr] \le C_{L} |\pi|.
\]
\end{Proposition}
\begin{pf}
Observe first that
\[
\mathbb{E}\bigl[\bigl|\tY^{\Re}_t - \tY^{\Re}_{\pi(t)}\bigr|^2 \bigr]
\le\mathbb
{E}\biggl[\biggl| \int_{\pi(t)}^t f(X_s,\tY^{\Re}_s,Z^{\Re}_s) \,\ud s +
\int_{\pi(t)}^t Z^{\Re}_s \,\ud W_s \biggr|^2 \biggr],\qquad 0\le t \le T.
\]
The proof is concluded combining this estimate with~(\ref{eqcontrolX}),
Propositions~\ref{estprioedsrdr}
and~\ref{prZbounded}.~%
\end{pf}
%

%As a direct consequence of the previous estimate, we deduce the
%following weaker

We now turn to the study of the regularity of the process $Z^{\Re}$.
\begin{Theorem}\label{ThmReguZ}
The process $Z^{\Re}$ satisfies
%
%e3.20 ###
%
\begin{equation}\label{EqReguZ1}
\mathbb{E}\biggl[\int_0^T | Z^{\Re}_s - \bar Z^{\Re}_{s} |^2 \,\ud s \biggr] \le C_L
(|\pi|^{1/2} + \kappa|\pi| ).\vadjust{\goodbreak}
\end{equation}
%
%point $a_0\in\Ic$ at time $0$, the optimal strategy $a$ associated to
%the discretely reflected BSDE $(Y,\tY,Z)$ satisfies
\end{Theorem}
%
% \romu{Il faut definir $\pi(\cdot)$ ici si on garde cet ordre dans le
%papier.}
%
\begin{pf}
A regularization argument as in proof of Proposition~\ref{prZbounded}
allows us to work under (Hr). From Remark~\ref{rebestapprox}, it is
clear that
%
%e3.21 ###
%
\begin{equation}\label{eqmajocontrol}
\mathbb{E}\biggl[\int_0^T | Z^{\Re}_s - \bar Z^{\Re}_{s} |^2 \,\ud s \biggr]
\le\mathbb{E}\biggl[\int_0^T \bigl| Z^{\Re}_s - Z^{\Re}_{\pi(s)} \bigr|^2
\,\ud s \biggr].
\end{equation}
For $s\le T$ and $a = (\alpha_{k}, \theta_{k})_{k\ge0} \in\cA^{\Re
}_{s,\ell}$, $\ell\in\cI$, we define $(V^{a}_{s,t})_{s\le t \le T}$ by
\begin{eqnarray*}
V^{a}_{s,t} &:=& \E_{t} \Biggl[ \nabla_{x}g^{a_T}(X_T) \Lambda^{a}_{s,T}
D_sX_T + \int_s^T \nabla_{x} f^{a_u}(X_u,\tY^{\Re}_u,Z^{\Re}_u)
\Lambda_{s,u}^{a} D_sX_u \,\ud u \\
&&\hspace*{139pt}{} - \sum_{k=1}^{N^{a}} \nabla_{x} c_{ \alpha_{j-1} , \alpha_j }
(X_{\theta_{k}}) \Lambda^{a}_{s,\theta_k}(D_sX)_{ \theta_k} \Biggr].
\end{eqnarray*}
%
%recalling $a = (\alpha_{k}, \theta_{k})_{k}$.
%
%%% pause commande
%%%ca va mieux on y retourne !
%
We now fix $\ell\in\cI$ and denote, for $u \le T$, by $ a^{u}
\in
\cA^{\Re}_{u,\ell}$ the optimal strategy associated to the
representation of $(\tY^{\Re}_{u})^{\ell}$, recalling (ii) in Theorem
\ref{THswR}.

Observe that, by definition, we have
%
%e3.22 ###
%
\begin{equation}\label{eqstrategy}
N^{ a^{t}} = N^{ a^{u}} \quad\mbox{and}\quad a^{t} = a^{u},\qquad
r_j \le t \le u < r_{j+1} , j< \kappa.
\end{equation}
Fix $i<n$, and deduce from Proposition~\ref{PropRepresZ} and (\ref
{eqstrategy}) that
%
%e3.23 ###
%
\begin{eqnarray}\label{eqcontrolZ1}
\mathbb{E}[|(Z^{\Re}_t)^\ell- (Z^{\Re}_{t_i})^\ell|^2] &=& \mathbb
{E}[|V^{ a^{t}}_{t,t} - V^{ a^{\ti}}_{t_i,t_i}|^2] \nonumber\\[-8pt]\\[-8pt]
&\le&2
( \mathbb{E}[|V^{ a^{\ti}}_{t,t} - V^{ a^{\ti}}_{t_i,t}|^2] +
\mathbb{E}[|V^{ a^{\ti}}_{t_i,t} - V^{ a^{\ti
}}_{t_i,t_i}|^2])\nonumber
\end{eqnarray}
%
%here
for $t\in[\ti,\tip)$. Combining (Hr),~(\ref{eqmajoDX}),
(\ref{ControlReguDsXt}),~(\ref{eqlambdasup}),~(\ref{eqlambdacontrol1})
and Cauchy--Schwarz inequality with the definition of $V^a$, we deduce
%
%e3.24 ###
%
\begin{equation} \label{eqcontrolZ2bis}
\mathbb{E}[|V^{ a^{\ti}}_{t,t} - V^{ a^{\ti}}_{t_i,t}|^2] \le
C_{L} |
\pi|^{1/2} ,\qquad \ti\le t \le\tip, i\le n .
\end{equation}

Since $V^{ a^{\ti}}_{t_i, .}$ is a martingale on $[\ti,\tip]$,
we obtain
%
%e3.25 ###
%
\begin{eqnarray}\label{eqcontrolZ3}
\mathbb{E}[|V^{ a^{\ti}}_{\ti,t} - V^{ a^{\ti}}_{\ti,\ti}|^2]
&\le&
\mathbb{E}[|V^{ a^{\ti}}_{\ti,\tip} - V^{ a^{\ti}}_{\ti,\ti}|^2]
\nonumber\\
&\le&
\mathbb{E}[|V^{ a^{\ti}}_{\tip,\tip}|^2 - |V^{ a^{\ti}}_{\ti
,\ti}|^2]
+\mathbb{E}[|V^{ a^{\ti}}_{\ti,\tip}|^2 - |V^{ a^{\ti}}_{\tip
,\tip}|^2]
\\
&\le& \mathbb{E}[|V^{ a^{\ti}}_{\tip,\tip}|^2 - |V^{ a^{\ti
}}_{\ti,\ti}|^2]
+ C_{L} |\pi|^{1/2} ,\qquad \ti\le t\le\tip,\nonumber
\end{eqnarray}
where the last inequality follows from~(\ref{eqcontrolZ2bis}).
Combining~(\ref{eqcontrolZ1}),~(\ref{eqcontrolZ2bis}), (\ref
{eqcontrolZ3}) and summing up over $i$, we obtain
\begin{eqnarray*}
&&\mathbb{E}\biggl[\int_{0}^{T} \bigl|(Z^{\Re}_t)^\ell- \bigl(Z^{\Re}_{ \pi
(t)}\bigr)^\ell\bigr|^2 \,\ud t\biggr] \\
&&\qquad\le C_{L}|\pi|^{1/2} + |\pi| \Biggl(
\mathbb{E}[|V^{ a^{r_{\kappa-1}}}_{T,T}|^2 - |V^{
a^{0}}_{0,0}|^2 ] + \sum_{j=1}^{\kappa-1}
(|V^{ a^{\rjm}}_{\rj,\rj}|^2 - |V^{ a^{\rj}}_{\rj,\rj}|^2)
\Biggr).
\end{eqnarray*}
Combined with~(\ref{eqmajoDX}) and~(\ref{eqlambdasup}), this
concludes the proof since $\ell$ is arbitrary.
\end{pf}

%stratgie optimale. Il
%devrait y avoir telescopage dans l equation precedente. Le terme en $
%que l'on force a recommencer le processus $a$ en $l$ chaque grid
%point. Vous confirmez? qu'on devrait pouvoir recuperer:}
% \[\label{EqReguZ2}
% \E[ \int_0^T | Z^{a_s}_s - Z^{a_{s}}_{\pi(s)} |^2 ] \le
%C_L |\pi|^{\frac{1}{2}} ou
% \E[ \int_0^T | Z^{a_s}_s - Z^{a_{\pi(s)}}_{\pi(s)} |^2 ] &
% \]

%s3.5 ###
\subsection{Extension}
\label{sseeuler}

%Recall that $\pi$ $:=$ $\{t_0=0,\ldots,t_n=T\}$ is a grid on the
%time interval $[0,T]$, such that $\Re\subset\pi$ and $|\pi|n\leq
%L$.

We shall approximate the process~$X$ by its Euler\break scheme~$X^\pi$, with dynamics
%
%e3.26 ###
%
\begin{equation}\label{eqEulerX}\quad
\Xp_t = X_{0} + \int_0^t b\bigl(X^\pi_{\pi(s)}\bigr) \,\ud s + \int_0^t
\sigma
\bigl(X^\pi_{\pi(s)}\bigr) \,\ud W_s,\qquad 0\le t\le T.
\end{equation}

Classically, we have the following upper-bound, uniformly in $\pi$:
%
%e3.27 ###
%
\begin{equation}\label{UBXpi}
\mathbb{E}\Bigl[\sup_{0\le t\le T} |\Xp_t | ^p\Bigr]^{1/p} \le
C_{L}^p ,\qquad p\geq2.
\end{equation}

The control of the error between $X$ and its Euler scheme $\Xp$ is
well understood (see, e.g.,~\cite{klopla92}) and we have
%
%e3.28 ###
%
\begin{equation}\label{eqErrorEuler}
\mathbb{E}\Bigl[\sup_{0\le t\le T} | X_t - \Xp_t | ^p\Bigr]^{1/p}
\le
C_{L}^p | \pi|^{1/2} ,\qquad p\geq2.
\end{equation}
In this context, we denote by $(Y^{eu},\tY^{eu},Z^{eu})$ the unique
solution of
the reflected BSDE $\cD(\Re, c(\Xp), f(\Xp,\cdot), g(\Xp))$. Our main
result here is the counterpart of Proposition~\ref{prregY} and
Theorem~\ref{ThmReguZ} when $X$ is replaced by $\Xp$.

\begin{Proposition}\label{prregYeZe}
The following hold:
\[
\mathbb{E}\biggl[\int_0^T\bigl|\tY^{eu}_t - \tY^{eu}_{\pi( t)}\bigr|^2 \,\ud t\biggr]
\le C_L
|\pi|
\]
and
\[
\mathbb{E}\biggl[\int_0^T | Z^{eu}_s - \bar
Z^{eu}_{s} |^2 \,\ud s \biggr] \le C_L (|\pi|^{1/2} + \kappa|\pi
| ) .
\]
\end{Proposition}
\begin{pf}
We only sketch the main step of the proof since it follows exactly the
same arguments as the one used to obtain Proposition
\ref{prregY} and Theorem~\ref{ThmReguZ}.\vspace*{8pt}

\textit{Step} 1. We use a kernel regularization argument which allows us
to work
under~(Hr). In this case, we observe that $\Xp$ belongs to
$\Lc_a^{1,2}$ and satisfies
\begin{eqnarray*}
D_s \Xp_t &=& \sigma\bigl(\Xp_{\pi(s)}\bigr) +
\int_s^t \nabla_x b\bigl(\Xp_{\pi(r)}\bigr)
D_s \Xp_{\pi(r)} \,\ud r\\
&&{} + \int_s^t \sum_{j=1}^q \nabla_x \sigma^j
\bigl(\Xp
_{\pi(r)}\bigr) D_s \Xp_{\pi(r)} \,\ud W^j_r
\end{eqnarray*}
for $s \le t$. One then checks (see~\cite{boucha08}, Remark 5.2,
for details) that
%
%e3.29 ###
%
\begin{eqnarray}
\label{eqmajoDXp1}
\Bigl\| \sup_{s\le T} | D_s \Xp| \Bigr\|_{_{\Sc^p}} & < & \infty,
\nonumber\\
\label{eqmajoDXp2}
{\sup_{s\le u }} \| D_s\Xp_t - D_s \Xp_u \|_{_{\mathbf{L}^p}} +
\Bigl\|\sup_{t \le s \le T } | D_t \Xp_s - D_u \Xp_s |
\Bigr\|_{_{\mathbf{L}^p}}
&\le& C_L^p | t-u |^{1/2} ,\hspace*{-35pt}\\
\eqntext{0\le u\le t\le T.}
\end{eqnarray}
It is also straightforward that $(Y^{eu},\tY^{eu},Z^{eu})$ is Malliavin
differentiable and satisfies~(\ref{EqDY}) with $\Xp$ instead of
$X$.\vspace*{8pt}

\textit{Step} 2. In order to retrieve the results of the proposition, one
then follows
exactly the same steps and arguments as the ones used in the
previous Sections~\ref{sserepz} and~\ref{sseregz}.
\end{pf}

%%%%%%%%%%%%%%%%%%%%%%
%%%%%%%%%%%%%%%%%%%%%%
%%%%%%%%%%%%%%%%%%%%%%
%%%%%%%%%%%%%%%%%%%%%%

%!TEX root = main.tex
%s4 ###
\section{A discrete-time approximation for discretely reflected
BSDEs}\label{SubSecScheme}

We present here a discrete time scheme for the approximation of the
solution of the discretely obliquely reflected BSDE~(\ref{eqdeBSDEDOR1}).

Recall that $\pi:=\{t_0=0,\ldots,t_n=T\}$ is a grid on the time
interval $[0,T]$,
%with modulus $|\pi|$ ($|\pi|$ $:=$ $\max_{0\leq i\leq
%n-1}|t_{i+1}-t_{i}|$),
such that $\Re\subset\pi$ and $|\pi|n\leq L$.
%including the reflection points of $\Re$.
In the sequel, the process $X$ is approximated by its Euler scheme
$X^\pi$ (see Section~\ref{sseeuler} for details).

%whose dynamics are given by
% \be\label{eqEulerX}
% \Xp_t &=& X_{0} + \int_0^t b(X^\pi_{\pi(s)}) \,\ud s + \int_0^t
% \ee
%
%Classically, we have the following uniform in $|\pi|$ upper-bound:
% , p\geq2 .
% \ee
% The control of the error between $X$ and its Euler scheme $\Xp$ is
%well understood, see e.g.~\cite{klopla92}, and we have
% \be\label{eqErrorEuler}
% \esp{ \sup_{0\le t\le T} | X_t - \Xp_t | ^p}^{1/p} &\le&
%C_{L}^p | \pi|^{1/2} , p\geq2 .
% \ee

%s4.1 ###
\subsection{A Euler scheme for discretely obliquely reflected BSDEs}

We introduce a Euler-type approximation scheme for the
discretely reflected BSDEs.

%Starting from the terminal condition $g(\Xp_T)$, classical discrete
%time schemes for BSDEs compute reversely in time an approximation of
%$Y$ on the grid points of $\pi$.
%These algorithms, described for example in the recent survey
%follow backward moves trying to mimic the dynamics of $Y$ in a forward
%perspective.
%For this purpose, we propose to rename posthumously this type of
%approximation as the moonwalk scheme
%and study here its extension to the approximation of BSDE with oblique
%reflections.

Starting from the terminal condition
\[
\Yp_T = \tYp_T := g(\Xp_T) \in\cC(\Xp_T),
\]
we compute recursively, for $i \le n-1$,
%
%e4.1 ###
%
\begin{equation}\label{scheme}
\cases{
\bZp_{\ti} = (\tip-\ti)^{-1}\mathbb{E}[\Yp_\tip(W_\tip
-W_\ti)'\mid\mathcal{F}_{\ti}],
\vspace*{2pt}\cr
\tYp_{\ti} = \mathbb{E}[\Yp_\tip\mid\mathcal{F}_{\ti
}] + (\tip- \ti) f(\Xp_\ti,\tYp
_\ti,\bZp
_{\ti}),
\vspace*{2pt}\cr
\Yp_{\ti} = \tYp_\ti\mathbf{1}_{ \{ \ti\notin\Re\} } + \cP
(\Xp_\ti,\tYp
_\ti)\mathbf{1}_{ \{ \ti\in\Re\} } .}
\end{equation}

This kind of backward scheme has been already considered when no
reflection occurs (see, e.g.,~\cite{boutou04}) and in the reflected
case (see, e.g.,~\cite{boucha08,majzha05,cha09}).
See also~\cite{bouelitou09} for a recent survey on the subject.

Combining an induction argument with the Lispchitz-continuity
of $f$, $g$ and the projection operator, one easily checks that
the above processes are square integrable and that the
conditional expectations are well defined at each step of the
algorithm.
\begin{Remark}
(i) This so-called ``moonwalk'' algorithm is given by an
implicit formulation, and one should\vspace*{1pt} use a fixed point argument to
compute explicitly $\tYp$ at each grid point.
%As observed in
%scheme, replacing the second equation in \reff{scheme} by
% \[
% \tYp_{\ti} &=& \EF{ \Yp_\tip+ (\tip- \ti) f(\Xp_\ti,\tYp_\tip,
% \]
% }\end{Remark}
%
% \begin{Remark}
%{\rm

(ii) In the two-dimensional case, Hamad{\`e}ne and Jeanblanc
\cite{hamjea07} interpret $Y^1-Y^2$ as the solution of a doubly
reflected BSDE. It is worth noticing that the solution of the
corresponding discrete time scheme developed by~\cite{cha09} for the
approximation of doubly reflected BSDE exactly coincides with
$(Y^{\Re,\pi})^{1}-(Y^{\Re,\pi})^{2}$ derived here.
\end{Remark}

For later use, we introduce the piecewise continuous time scheme
associated to the triplet $(\Yp,\tYp, \bZp{})$. By the martingale
representation
theorem, there exists $\Zp\in{ \mathcal{H}^{2}}$ such that
\[
\Yp_\tip=\mathbb{E}_{\ti} [\Yp_\tip]+\int_\ti
^\tip\Zp_u \,\ud W_u ,\qquad
i\le n-1,
\]
and by the It\^o's isometry, for $i \le n-1$,
%
%e4.2 ###
%
\begin{equation}\label{eqrepbZp}
\bZp_\ti= \frac{1}{\tip-\ti} \mathbb{E}\biggl[\int_\ti^\tip\Zp
_s \,\ud s\Bigm|\mathcal{F}_{\ti}\biggr].
\end{equation}
We set $\bZp_t := \bZp_{\pi(t)}$ for $t\in[0,T]$, define $\tYp$ by
%
%e4.3 ###
%
\begin{eqnarray} \label{eqdeftYpcont}
\tYp_t &=& \Yp_\tip
+ (\tip-t) f(\Xp_\ti,\tYp_\ti,\bZp_\ti)\nonumber\\[-8pt]\\[-8pt]
&&{}- \int_t^\tip\Zp_u\,\ud W_u,\qquad \ti\le t\le\tip, i\in
\Ic,\nonumber
\end{eqnarray}
and introduce $\Yp$ on $[0,T]$ by
$
\Yp_t :=\tYp_t\mathbf{1}_{ \{ t \notin\Re\} } + \cP(\Xp
_{t},\tYp_t)\mathbf{1}_{ \{ t \in\Re\} } .
$

This can be rewritten as
%
%e4.4 ###
%
\begin{equation} \label{eqdeftYpcont2}\quad
\cases{\displaystyle \tYp_{t }
= g(\Xp_{T}) + \int_{t}^T f\bigl(\Xp_{\pi(u)},\tYp_{\pi(u)},\bZp
_{u}\bigr)\,\ud u -
\int_{t}^T \Zp_u \,\ud W_u\vspace*{2pt}\cr
\hphantom{\tYp_{t }
=}{} + \displaystyle (K ^{\Re,\pi}_T-K ^{\Re,\pi}_t),
\vspace*{2pt}\cr
\displaystyle K ^{\Re,\pi}_t := \sum_{r \in\Re\setminus{\{ 0 \}}} \Delta K
^{\Re
,\pi}_r \mathbf{1}_{ \{ r \le t \} } \quad\mbox{and }\vspace*{2pt}\cr
\Delta K ^{\Re,\pi}_t:= \Yp_t -
\tYp_t
= -(\tYp_t - \tYp_{t-}),
\vspace*{2pt}\cr
\Yp_{t} = \tYp_t\mathbf{1}_{ \{ t \notin\Re\} } + \cP(\Xp
_{t},\tYp_t)\mathbf{1}_{ \{ t \in\Re\} } , \qquad\mbox{$0\le t \le
T$}.}
\end{equation}

We finally provide a useful a priori estimate for the solution of
the discrete time scheme whenever $f$ does not depend on $z$, whose
proof is postponed until Appendix~\ref{secap2}.
\begin{Proposition}\label{4BUShceme}
If $f$ does not depend on $z$ and $|\pi|L<1$, the following
bound holds:
%
%e4.5 ###
%
\begin{equation}\label{4propUBScheme}
\E\Bigl[ \sup_{0\leq i\leq n} |\tYp_{t_{i}}|^p\Bigr] \leq
C_{L}^p,\qquad
p\geq2.
\end{equation}
Recall that $C_{L}^p$ neither depends on $\Re$ nor on $\pi$.
\end{Proposition}

%s4.2 ###
\subsection{Convergence results}
The next proposition provides a control on the error between the
discrete-time scheme~(\ref{scheme}) and the solution of the
discretely reflected BSDE~(\ref{eqdeBSDEDOR1}).

% the discretely reflected BSDE in the general setting.
%
\begin{Proposition}\label{ThmErrorScheme}
The following holds: %discrete-time scheme $(\tYp,\bZp)$ converges to $(
%
%e4.6 ###
%
\begin{eqnarray}\label{eqThmErrorScheme}\quad
&&\sup_{t\in[0,T]} \mathbb{E}[|\tY^{\Re}_t - \tYp_t|^2 +
|Y^{\Re}_t - \Yp_t|^2] + \mathbb{E}\biggl[\int_0^T|Z^{\Re}_s - \bZp
_{s} |^2 \,\ud s\biggr]\nonumber\\[-8pt]\\[-8pt]
&&\qquad\le
C_{L} |L_\cP|^{2\kappa}( |\pi|^{1/2}+\kappa|\pi|),\nonumber
\end{eqnarray}
where we recall that $L_\cP=\sqrt{d}$ is the Lipschitz constant of the
projection operator~$\cP$.
\end{Proposition}
\begin{pf}
As in Section~\ref{sseeuler}, we consider
$(Y^{eu},\tY^{eu},Z^{eu})$ the unique solution of
the reflected BSDE $\cD(\Re, c(\Xp), f(\Xp,\cdot), g(\Xp))$.
Using Proposition~\ref{prstabY}, the Lipschitz property of $f$,
$g$, $c$ and~(\ref{eqErrorEuler}), we obtain
%
%e4.7 ###
%
\begin{eqnarray}\label{eqpourromu1}
&&\sup_{t\in[0,T]} \mathbb{E}[|\tY^{\Re}_t - \tY^{eu}_t|^2 +
|Y^{\Re}_t - Y^{eu}_t|^2] + \frac1\kappa
\mathbb{E}\biggl[\int_0^T|Z^{\Re}_s - Z^{eu}_{s} |^2 \,\ud s\biggr]\nonumber\\[-8pt]\\[-8pt]
&&\qquad \le C_L |\pi|.\nonumber
\end{eqnarray}

Using the same arguments as in the proof of Proposition 3.4.1,
Step 1.a in~\cite{cha08}, for example, we get the following inequality:
%
%e4.8 ###
%
\begin{eqnarray}\label{eqcontrolschemegene}
&&\sup_{t\in[\ti,\tip)} \mathbb{E}[|\tY^{eu}_t - \tYp_t|^2 +
|Y^{eu}_t - \Yp_t|^2] \nonumber\\
&&\quad{}+ \mathbb{E}\biggl[\int_\ti^\tip|Z^{eu}_s - \bZp
_{s} |^2 \,\ud s\biggr] \nonumber\\[-8pt]\\[-8pt]
&&\qquad \le
C_L \biggl( \mathbb{E}\biggl[|Y^{eu}_\tip- \Yp_\tip|^2\nonumber\\
&&\qquad\hphantom{\le
C_L \biggl( \mathbb{E}\biggl[}
{} + \int_\ti^\tip
\bigl(\bigl|\tY^{eu}_s - \tY^{\Re}_{\pi(s)}\bigr|^2 +\bigl|Z^{eu}_s - \bar Z^{\Re
}_{\pi(s)}\bigr|^2\bigr) \,\ud s\biggr]
\biggr).\nonumber
\end{eqnarray}

There\vspace*{1pt} are two differences with the proof of Proposition 3.4.1 in
\cite{cha08}. First, $\cP$~here depends both on $x$ and $y$: but
this is not a problem since $(Y^{eu},\tY^{eu},\break Z^{eu})$ and $(\Yp
,\tYp,
\Zp)$
are parametrized by the same forward process~$\Xp$.

Second, $\cP$ is not $1$-Lipschitz but only $L_\cP$-Lipschitz, with
$L_\cP
> 1$, in its~$y$ component. This explains the term $|L_\cP|^{2\kappa}$
in~(\ref{eqThmErrorScheme}). Indeed, we have, for $i < n$,
\[
|Y^{eu}_\tip- \Yp_\tip|^2 = |\cP(\Xp_\tip,\tY^{eu}_\tip) -
\cP(\Xp_\tip,\tYp_\tip) |^2 \le|L_\cP|^2 |\tY^{eu}_\tip- \tYp
_\tip|^2.
\]
This leads, using an induction argument (see, e.g., Step 1.b in the
proof of Proposition 3.4.1 in~\cite{cha08}), to
\begin{eqnarray*}
&&\sup_{t\in[0,T]} \mathbb{E}[|\tY^{eu}_t - \tYp_t|^2 + |Y^{eu}_t -
\Yp_t|^2 ]
+ \mathbb{E}\biggl[\int_0^T |Z^{eu}_s - \bZp_{s} |^2 \,\ud s \biggr]\\
&&\qquad \le
C_L |L_\cP|^{2 \kappa} \biggl( |\pi| + \int_0^T\bigl(\bigl|\tY^{eu}_s - \tY
^{\Re
}_{\pi(s)}\bigr|^2 +\bigl|Z^{eu}_s - \bar Z^{\Re}_{\pi(s)}\bigr|^2\bigr)
\,\ud s \biggr).
\end{eqnarray*}

Combining the last inequality with Proposition~\ref{prregYeZe} and
(\ref{eqpourromu1}) completes the proof.\vspace*{-3pt}
\end{pf}

The term $|L_\cP|^{2\kappa}$, even when $\kappa$ is small can be very
large. Moreover, we shall see in the next section that it yields to
a poor convergence rate for continuously reflected BSDEs. This term
is due to the ``geometric'' approach, used in the proof of
Proposition~\ref{ThmErrorScheme}, and the fact that $\cP$ is only
$L_\cP$-Lipschitz with $L_\cP>1$. We obtain below a better control, using
the stability results proved at the end of Section~\ref{sedisc} but
unfortunately
under the assumption that~$f$ does not depend on $z$. The optimal
choice for $\kappa$ in terms of $|\pi|$ is discussed in Section \ref
{SubSecConvergence} below.\vspace*{-3pt}
\begin{Theorem}\label{THcvSH}
If $f$ does not depend on $z$, the following holds:
\begin{eqnarray*}
\sup_{t\in[0,T]} \mathbb{E}[|\tY^{\Re}_{t}-\tYp_{t}|^2 + |Y^{\Re
}_{t}-\Yp_{t}|^2] & \le & C_L|\pi|,\\[-2pt]
\mathbb{E}\biggl[\int_{0}^{T}|Z^{\Re}_{t} - \bZp_{t}|^2 \,\ud t\biggr] & \le &
C_{L}(\kappa|\pi| + |\pi|^{1/2} )
\end{eqnarray*}
for $|\pi|$ small enough.\vspace*{-3pt}
% \jef{distinguer avec Hf ajouter ici les resultats sur le Z}
\end{Theorem}
\begin{pf}
We use here the stability results of Proposition~\ref{prstabY}
setting $({}^{1}Y^{\Re},{}^{1}\tY^{\Re},{}^{1} Z^{\Re}) = (Y^{\Re
}, \tY
^{\Re}, Z^{\Re})$ with ${}^{1}F\dvtx(s,y,z) \mapsto f(X_{s},\tY^{\Re}_{s})$
and\vspace*{1pt} $({}^{2}Y^{\Re},\break{}^{2}\tY^{\Re}, {}^{2}Z^{\Re}) = (\Yp, \tYp,
\Zp
)$, with ${}^{2}F\dvtx(s,y,z) \mapsto f(\Xp_{\pi(s)},\tYp_{\pi(s)})$.
Combining\vspace*{1pt}~(\ref{4propUBScheme}) and Proposition~\ref{estprioedsrdr}
with the Lipschitz property of $f$, it is clear that (C$_4$)
holds. Applying Proposition~\ref{prstabY} and~(\ref{eqErrorEuler}),
we derive, for $t\in[0,T]$,
%
%e4.9 ###
%
\begin{eqnarray}\label{eqconvscheme1}\qquad
&&
\E|\tY^{\Re}_{t}-\tYp_{t}|^2 + \frac{1}{\kappa} \int_t^T\E
|Z^{\Re
}_{s}-\Zp_{s}|^2\,\ud s\nonumber\\[-9pt]\\[-9pt]
&&\qquad\le
C_{L}\biggl(|\pi|+ \int_{t}^T\E\bigl|\tY^{\Re}_{s} -\tY^{\Re}_{{\pi
(s)}}\bigr|^2
\,\ud s+\int_{t}^T\E\bigl|\tYp_{\pi(s)}-\tY^{\Re}_{{\pi(s)}}\bigr|^2\,\ud s\biggr).\nonumber
\end{eqnarray}
Applying the discrete version of Gronwall's lemma to estimate (\ref
{eqconvscheme1}) rewritten at time $t=t_j\in\pi$, we deduce
%The same reasoning applied at time $t_{j}\in\pi$ with $t_{j}\geq t$
%leads to
%+ \int_{t}^T\E|\tY_{s}-\tY_{{\pi(s)}}|^2ds+ |\pi| \sum_{k=j+1}^n \E|
% Applying the discrete version of Gronwall's lemma, we deduce
%
%e4.10 ###
%
\begin{eqnarray}\label{eqconvscheme1bis}
\E|\tY^{\Re}_{t_{j}}-\tYp_{t_{j}}|^2 \le C_{L}\biggl(|\pi|
+ \int_{t}^T\E\bigl|\tY^{\Re}_{s}-\tY^{\Re}_{{\pi(s)}}\bigr|^2\,\ud s \biggr)
,\nonumber\\[-9pt]\\[-9pt]
\eqntext{0\le t\le t_{j}\le T, t_j\in\pi.}
\end{eqnarray}
Plugging this estimate into~(\ref{eqconvscheme1}), we compute
\begin{eqnarray*}
&&\E|\tY^{\Re}_{t}-\tYp_{t}|^2 + \frac{1}{\kappa} \int_t^T
\E
|Z^{\Re}_{s}-\Zp_{s}|^2\,\ud s\\
&&\qquad\le
C_{L}\biggl(|\pi| + \int_{t}^T\E\bigl|\tY^{\Re}_{s}-\tY^{\Re}_{{\pi
(s)}}\bigr|^2\,\ud s\biggr) ,\qquad 0\le t\le T,
\end{eqnarray*}
which combined with Proposition~\ref{prregY} leads to the first claim
of the theorem.

Observe from the representations~(\ref{defbarZ}) and~(\ref{eqrepbZp}) that
\begin{eqnarray*}
&&\mathbb{E}\biggl[\int_{0}^{T} |Z^{\Re}_{t} - \bZp_{t}|^{2} \,\ud t\biggr]
\\
&&\qquad\le
C_{L} \biggl(\mathbb{E}\biggl[\int_{0}^{T} |Z^{\Re}_{t} - \bar{Z}^{\Re
}_{t}|^{2} \,\ud t\biggr] +\mathbb{E}\biggl[\int_{0}^{T} |Z^{\Re}_{t} - \Zp
_{t}|^{2} \,\ud t\biggr] \biggr).
\end{eqnarray*}
Plugging~(\ref{EqReguZ1}), estimate~(\ref{eqconvscheme1}) written
at time $t=0$ and the first claim of this theorem into this expression
concludes the proof.
%Using Proposition~\ref{prstabY}, Proposition~\ref{prregY},
%that
% .
% The proof of the second claim of the Theorem is obtained combining
%the last inequality with \eqref{eqconvscheme2} and Theorem
\end{pf}

%%%%%%%%%%%%%%%%%%%%%%
%%%%%%%%%%%%%%%%%%%%%%
%%\input{continuous}
%%%%%%%%%%%%%%%%%%%%%%
%%%%%%%%%%%%%%%%%%%%%%

%!TEX root = main.tex

%s5 ###
\section{Extension to the continuously reflected case}
\label{secontinuous}

In this section, we extend the convergence results of the scheme (\ref
{scheme}) to the case of continuously reflected BSDEs. To this end, we
show that the error between discretely and continuously obliquely
reflected BSDEs is controlled in a convenient way.

%s5.1 ###
\subsection{Continuously obliquely reflected BSDEs}

In the sequel, we shall use the following assumption on $f$:

\begin{longlist}[(ii) (Hz)]
\item[(i) (Hz)]
The function $f$ is bounded in its last variable:
${\sup_{z
\in\cM^{d,q}}} |f(0,0,\break z)| \leq C_L$
and the following assumption on the cost $c$.

\item[(ii) (Hc)] For\vspace*{1pt} $i,j \in\cI$, the function $c^{ij}$ is equal
to $^{1}c^{ij} - ^{2}c^{ij}$, with $^{1}c^{ij}$ is~$C^{2}$ with bounded
first and second derivatives and $^{2}c^{ij}$ is a convex function with
bounded first derivative.
\end{longlist}

This last assumption is needed to retrieve some regularity
on the reflecting process $K$ (see Lemma~\ref{lereguK} below).

% We now study the behavior of the solution of the discretely reflected
% BSDE as the time mesh of the reflection grid goes to zero.
% We prove that it converges to the solution of the corresponding
% continuously reflected BSDE, and provide a control on the quadratic
% distance between both of them.

%For this purpose, we introduce an increasing sequence of grids
%$(\Re^n)_{n\in\N}$, partitions of $[0,T]$, such that $\cup_n \Re^n$ is
%dense in $[0,T]$. The time mesh $|\Re^n|$ associated to each grid
%$\Re^n$ is the biggest gap between two consecutive grid points.
%We denote by $(\tY^n,Y^n,Z^n,K ^n)$ the solution to the associated
%discretely reflected BSDE, and by
We denote by
$(Y ,Z ,K )\in({\Sc^2}\times{\Hc^2}\times\mathbf{A}^2)^\Ic$
the solution of the continuously obliquely reflected BSDE $ \cC
([0,T],c(X),f(X,\cdot), g(X_{T}))$ defined by
%
%e5.1 ###
%
\begin{equation}\label{eqBSDECOR}\qquad
\cases{
\displaystyle Y ^i_t = g^i(X_T)+\int_t^T f^i(X_s,Y _s^i, Z _s^i)\,\ud s-\int
_t^T Z
^i_s \,\ud W_s + K ^i_T - K ^i_t,
\vspace*{2pt}\cr
\displaystyle Y ^i_t \geq\max_{j\in\Ic} \{Y ^j_t - c^{ij}(X_{t})\} ,
\qquad\hspace*{64.6pt}
\mbox{$0\leq
t\leq T$},
\vspace*{2pt}\cr
\displaystyle \int_0^T\Bigl[Y _t^i-\max_{j\in\Ic}\{Y ^j_t - c^{ij}(X_{t})\}\Bigr] \,\ud
K
^i_t=0, \qquad
\mbox{$i\in\Ic$}.}
\end{equation}
Under the assumption on $f$, $g$ and $c$, the existence and
uniqueness of such a solution is given in~\cite{hz08,hutan08}.

The solution of~(\ref{eqBSDECOR}) has also a representation property
in term of switched BSDEs, recalling~(\ref{eqUa}). Here, of course,
the switching times of the strategy are not restricted to take their
values in $\Re$. We refer to~\cite{chaelikha10} for more details.
\begin{Theorem}\label{thsummary}
There exists, for any fixed initial condition
$(t,i)\in\break[0,T]\times\Ic$, an optimal switching strategy $\dot
a:=(\dot\theta_k,\dot\alpha_k)_{k\ge0}\in\cA_{t,i}$, such that
%
%e5.2 ###
%
\begin{equation}\label{Snellcontinuous}
Y ^{i}_{t} = U^{\dot a}_{t} = \mathop{\esssup}_{a\in\Ac_{t,i}} U^a_t,
\qquad\mathbb{P}\mbox{-a.s.}
\end{equation}
\end{Theorem}

We deduce from~(\ref{Snellcontinuous}), Theorem~\ref{THswR}(iii), the
monotonicity property of $\cP$ and
(\ref{eqBSDECOR})
%
%e5.3 ###
%
\begin{equation}\label{eqordering}
Y \succeq Y^{\Re} \succeq\tY^{\Re}\qquad \mbox{for any grid }
\Re.
\end{equation}

Moreover, most of the estimates presented in Section \ref
{sedisc} for discretely reflected BSDEs hold true for continuously
reflected BSDEs. For reader's convenience, we collect them in the
following proposition. The proof itself is postponed to Appendix
\ref{secap3}.
\begin{Proposition} \label{prstabcont} The following a priori
estimates hold.
For any $p \ge2$,
%
%e5.4 ###
%
\begin{eqnarray}\label{eqestimateia}
|Y _{t}|^p + \mathbb{E}_{t} \biggl[\biggl(\int_{t}^T |Z _{s}|^2
\,\ud s \biggr)^{p/2} \biggr] + \mathbb{E}_{t} [|K
_{T}-K _{t}|^p] \le\mathbb{E}_{t} [\beta^X] ,\nonumber\\[-8pt]\\[-8pt]
\eqntext{0 \le t \le T,}%\|Y \|_{_{{\Sc^p}}} + \|Z \|_{_{{\Hc^p}}} + \|
%K _T\|_{_{{\bf L^p}}}
\end{eqnarray}
and, for all $(t,i) \in[0,T]\times\cI$, the optimal strategy $\dot{a}
\in\cA_{t,i}$ satisfies
%, for all $(t,i) \in[0,T]\times\cI$,
%
%e5.5 ###
%
\begin{equation} \label{eqestimateib}
\mathbb{E}_{t} \Bigl[\sup_{s \in[t,T]} |U^{\dot a}_{s}|^p
\Bigr] + \mathbb{E}_{t} [|N^{\dot a}|^p] \le\mathbb
{E}_{t} [\beta^{X}].
\end{equation}
%
% INUTILE MAINTENANT
%
%Z , {}^\ellK)$ be the solution of the continuously obliquely
%reflected BSDE $ \cC([0,T],c({}^\ell X),f({}^\ell X,\cdot), {}^\ell
%that ${}^{1}\xi\vee{}^2\xi+ \sup_{s \in[0,T]} |{}^1X_{s}|\vee
%|{}^2X_{s}| \le\beta^X $, then we have
%|{}^1 Y _{t} - {}^2 Y _{t}|^2 \le C_{L} \EFp{t}{|{}^1 \xi-{}^2
%X_{s} - {}^2 X_{s}|^4}^{1/2}, 0 \le t \le T .
\end{Proposition}
\subsection{Error between discretely and continuously reflected BSDEs}

We first provide a control of the error on the grid points of $\Re$
between the solutions of the obliquely discretely and continuously
reflected BSDEs~(\ref{eqBSDEDOR}) and~(\ref{eqBSDECOR}).

\begin{Theorem}\label{ThmErrorDiscretize}
Under \textup{(Hz)}, the following holds: %for $p\ge2$
%and $\varepsilon> 0$
%and under $\HYP{f}$
%
%e5.6 ###
%
\begin{equation}\label{4errordisccontcasgen}
\mathbb{E}\Bigl[\sup_{r \in\Re} \{|Y _r - \tY^{\Re}_r|^2+|Y _r -
Y^{\Re}_r|^2\}\Bigr] \le C_L^\eps|\Re|^{1-\eps}
,\qquad \eps>0.
\end{equation}
%
%and for all $\eps>0$, there exits a constant $C_{L}^\eps$ such that
% .
Moreover, if the cost functions are constant, the last inequality holds
true with $\varepsilon= 0$.
% \sup_{r \in\Re}\esp{ \{|Y _r - \tY_r|^2+|Y _r - Y_r|^2
% \]
% and
% \be\label{4errordisccontcaspart}
% \esp{ \sup_{r \in\Re}\{|Y _r - \tY_r|^2+|Y _r - Y_r|^2
% .
\end{Theorem}
%
%[FAIRE UNE REMARQUE SUR LA Convergence(?) sans vitesse si pas HF ??]
%
%We already know from \reff{eq tYn tY bY} that $\tY^i \leY^i$ $\Pas$
%for $1\le i\le d$. Let fix $(t,i)\in[0,T]\times\{1,\ldots,d\}$ and
%try to prove the opposite inequality.
%
%We naturally denote by $\bA$ and $(A^n)_n$ the cumulative costs
%associated respectively to the switching strategy $\ba$ and $(a^n)_n$.
%We provide in the next lemma a useful estimate on these processes.
%
\begin{pf}
The proof of this result relies mainly on the interpretation in
terms of switched BSDEs provided in Section~\ref{subsediscopt}. For a
fixed $(t,i)\in{[0,T]\times\Ic}$, we associate to the optimal strategy
$\dot a=(\dot\theta_{k},\dot\alpha_{k})_{k}\in\Ac_{t,i}$ not
restricted to lie in the grid $\Re$, the corresponding ``discretized''
strategy $a :=(\theta_k,\alpha_k)_{k\ge0}\in\Ac^{\Re}_{t,i}$
defined by
%
%e5.7 ###
%
\begin{equation}\label{defthetaalpha}
\theta_k := \inf\{ r\ge\dot\theta_k ; r\in\Re
\} \quad\mbox{and}\quad \alpha_k := \dot{\alpha}_k ,\qquad k\ge0 .
\end{equation}

\vspace*{4pt}

\textit{Step} 1. We first derive two key controls on the distance
between $A^{\dot a}$ and~$A^{a}$.\vadjust{\goodbreak}

We fix $p\geq2$ and, since $\dot\theta_{k}\le\theta_{k}$, $k\geq1$,
we compute
%
%e5.8 ###
%
\begin{eqnarray}\label{4decompdeltaA}
&&\biggl(\int_{t}^T | A^{\dot a}_s - A^a_s|^2 \,\ud s\biggr)^{p/2}\nonumber\\
&&\qquad =
\Biggl(\int_{t}^T \Biggl|
\sum_{k=1}^{ N^{\dot a}} c^{\dot\alpha_{k-1}\dot\alpha
_k}(X_{\dot\theta_{k}}){\mathbf1}_{\dot\theta_k \le s }
- c^{\dot\alpha_{k-1}\dot\alpha_k}(X_{\theta_{k}}){\mathbf
1}_{\theta_k \le s}
\Biggr|^2 \,\ud s\Biggr)^{ p/2} \nonumber\\[-8pt]\\[-8pt]
&&\qquad \le
C_{L}^p\int_{t}^T \Biggl|
\sum_{k=1}^{N^{\dot a}}[
c^{\dot\alpha_{k-1}\dot\alpha_k}(X_{\theta_{k}})-c^{\dot\alpha
_{k-1}\dot\alpha_k}(X_{\dot\theta_{k}})
]{\mathbf1}_{\theta_k \le s }
\Biggr|^p\,\ud s \nonumber\\
&&\qquad\quad{} +
C_{L}^p\Biggl(\int_{t}^T \Biggl|
\sum_{k=1}^{N^{\dot a}}
c^{\dot\alpha_{k-1}\dot\alpha_k}(X_{\dot\theta_{k}}){\mathbf
1}_{\dot
\theta_k \le s < \theta_{k} }
\Biggr|^2 \,\ud s\Biggr)^{ p/2}.\nonumber
\end{eqnarray}
%
%where we set $\dot N$ $:=$ $N^{\dot a}$. % denotes the almost finite
%number of switch in the strategy $\ba$.
%[changer en $\bar{N}$]
Using the convexity inequality $(\sum_{k=1}^n |x_{k}|)^p\leq
n^{p-1}\sum_{k=1}^n |x_{k}|^{p}$, we obtain
%
%e5.9 ###
%
\begin{eqnarray}\label{4maj1deltaA}
&&\Biggl(\int_{t}^T \Biggl|
\sum_{k=1}^{N^{\dot a}}
c^{\dot\alpha_{k-1}\dot\alpha_k}(X_{\dot\theta_{k}}){\mathbf
1}_{\dot
\theta_k \le s < \theta_{k} }
\Biggr|^2 \,\ud s\Biggr)^{ p/2} \nonumber\\[-8pt]\\[-8pt]
&&\qquad\leq C^p_{L} \Bigl(1+\sup_{t\in
[0,T]}|X_{t}|^p\Bigr)|N^{\dot a}|^p|\Re|^{p/2}.\nonumber
\end{eqnarray}
Using once again the same convexity inequality with $p=2$, the
Lipschitz property of the maps $(c^{ij})_{i,j\in\Ic}$ and the
definition of $\dot\theta_{k}$ and $\theta_{k}$, we get
\begin{eqnarray*}
\int_{t}^T \Biggl|
\sum_{k=1}^{N^{\dot a}}[
c^{\dot\alpha_{k-1}\dot\alpha_k}(X_{\theta_{k}})-c^{\dot\alpha
_{k-1}\dot\alpha_k}(X_{\dot\theta_{k}})
] {\mathbf1}_{\theta_k \le s }
\Biggr|^p \,\ud s
& \leq&
C^p_{L}|N^{\dot a}|^{p-1} \sum_{k=1}^{N^{\dot a}}|X_{\theta
_{k}}-X_{\dot\theta_{k}}|^p \\
& \leq&
C^p_{L}|N^{\dot a}|^{p} \chi^{|\Re|,p},
\end{eqnarray*}
where $\chi^{|\Re|,p} := \sum_{k=1}^{\kappa}\sup_{r\in
[r_{k-1},r_{k}]}|X_{r}-X_{r_{k}}|^p$.\vspace*{1pt}

Plugging this estimate and~(\ref{4maj1deltaA}) in (\ref
{4decompdeltaA}), we deduce
%
%e5.10 ###
%
\begin{eqnarray}\label{eqcontrolbA-A1}
&&\biggl(\int_{t}^T | A^{\dot a}_s -A^a_s |^2 \,\ud s
\biggr)^{p/2} \nonumber\\[-8pt]\\[-8pt]
&&\qquad\le
C^p_L |N^{\dot a}|^p\Bigl(\Bigl(1+\sup_{s\in[0,T]}|X_{s}|^p
\Bigr)|\Re
|^{p/2}+ \chi^{|\Re|,p}\Bigr).\nonumber
\end{eqnarray}

Observe also that, for $r \in\Re$, we have ${\mathbf1}_{\dot\theta
_{k}\leq
r}={\mathbf1}_{\theta_{k}\leq r}$ which gives
%
%e5.11 ###
%
\begin{eqnarray}\label{eqcontrolbA-A2}
|A^{\dot a}_r - A^a_r|^p &\le&\Biggl(\sum_{k=1}^{N^{\dot a}}
|c^{\dot\alpha_{k-1}\dot\alpha_k}(X_{\dot\theta_{k}})-c^{\dot
\alpha_{k-1}\dot\alpha_k}(X_{\theta_{k}})|{\mathbf
1}_{\theta_{k}\le r}
\Biggr)^p\nonumber\\[-8pt]\\[-8pt]
&\le& C_L |N^{\dot a}|^p \chi^{|\Re|,p}.\nonumber
\end{eqnarray}

\vspace*{4pt}

%
% CF CORO et preuve EN APPENDICE MAINTENANT
%switches $N^{\dot a}$.
%
%Recall that $(\dot Y^{\dot a},\dot Z^{\dot a})$ satisfies the switched
%equation \reff{eq Ua} with $a=\dot a$ starting at time $t$ in state
%$i$.
%Therefore, using the Lipschitz property of $f$ and Burkholder Davis
%Gundy inequality, we get
% \begin{eqnarray*}
% \E_{t}[|A_T^{\dot a}|^{2p} ] & \le& C_L^{p} \E_{t}[
%1 + \sup_{0\le s \le T} |X_s|^{2p} + \sup_{0\le s \le T} |\dot Y^{\dot
%a_s}_s|^{2p} + (\int_t^T | \dot Z_{s}^{\dot a_s}|^2 \,\ud s)^{p}
%+ |A^{\dot a}_t|^{2p} ] .
% \]
% But the cost structure \eqref{CostStructure} implies $|N^{\dot a}|
%and $|A^{\dot a}_t|\leq C_L(1+|X_t|)$ leads to
% \be\label{ControlNadot}
% \E_{t}[|N^{\dot a}|^{2p} ] & \le& C_L^{p} \E_{t}[
%1 + \sup_{0\le s \le T} | X_s|^{2p} + \sup_{0\le s \le T} |\dot
%Y_s|^{2p} + (\int_t^T | \dot Z_s|^2 ds)^{p} ] .
% \ee

\textit{Step} 2. We now prove the main result of the
theorem.

We introduce the processes $\Gamma:= U^{a} - A^{a}$ and $\dot\Gamma
:= U^{\dot a} - A^{\dot a}$, so that
% the
% dynamics of $\dot\Gamma-\Gamma$ on $[t,T]$ is given by
% \begin{eqnarray*}
% \dot\Gamma_t - \Gamma_t &= & \int_t^T [
% f^{\dot a_s}(X_{s},\dot\Gamma_s + A^{\dot a}_s, V^{\dot a}_s)
% - f^{a_s}(X_{s},\Gamma_s + A^a_s, V^{a}_s) ] \,\ud s - \int_t^T
%(V^{\dot a}_s - V^{a}_s) .\,\ud W_s .
% \]
% Observe that,
%
%e5.12 ###
%
\begin{equation}\label{eqteubee1}
|U^a - U^{\dot a} | \le|\Gamma- \dot\Gamma| + |A^a - A^{\dot a}|.
\end{equation}
%
%A^{n}_s, V^{a^n}_s)ds - \int_t^T V^{a^n}_s . dW_s , \\
Applying It\^o's formula to the continuous process $|\dot\Gamma-
\Gamma|^2$ on $[t,T]$, using Gronwall's lemma and the Lipschitz
property of $f$, we obtain
%
%e5.13 ###
%
\begin{eqnarray}\label{4eqpuissance2}
&&| \dot\Gamma_{t} - \Gamma_{t}|^2 \nonumber\hspace*{-35pt}\\[-9pt]\\[-9pt]
&&\qquad\le C_L \E_{t}\biggl[\int_{t}^T
|[f^{\dot a_s} - f^{a_s}](X_{s},U^{\dot a}_s, V^{\dot a}_s)
|^2 \,\ud s + \int_{t}^T | A^{\dot a}_s-A^a_s|^2 \,\ud s\biggr].\nonumber\hspace*{-35pt}
\end{eqnarray}
Elevating this expression to the power $p \over2$, we deduce
%
%e5.14 ###
%
\begin{eqnarray} \label{eqteubee2}
&&| \dot\Gamma_{t} - \Gamma_{t}|^p \nonumber\\[-2pt]
&&\qquad\le C^p_L \E_t\biggl[\biggl(\int
_{t}^T |[f^{\dot a_s} - f^{a_s}](X_{s},U^{\dot a}_s, V^{\dot
a}_s)|^2 \,\ud s\biggr)^{p/2}\\[-2pt]
&&\qquad\hspace*{100.6pt}{} +\biggl( \int_{t}^T |
A^{\dot
a}_s-A^a_s|^2 \,\ud s\biggr)^{p/2}\biggr].\nonumber
\end{eqnarray}
Combining the definition of $\theta$ with the Lipschitz property of $f$
and (Hz), we compute
\begin{eqnarray*}
&&\int_{t}^T |[f^{\dot a_s} -
f^{a_s}](X_{s},U^{\dot a}_s, V^{\dot a}_s)|^2 \,\ud s\\[-2pt]
&&\qquad =
\int_{t}^T \Biggl|\sum_{k=1}^{ N^{\dot a}} f^{\alpha
_{k-1}}(X_s,U^{\dot a}_s, V^{\dot a}_s)
( {\mathbf1}_{\dot\theta_{k-1} \le s < \dot\theta_k} -
{\mathbf1}_{\theta_{k-1} \le s < \theta_k} ) \Biggr|^2 \,\ud s\\[-2pt]
&&\qquad \le
C_L |N^{\dot a}|^2 \sup_{s \in[0,T]} (1+|X_s|^2 + | U^{\dot
a}_s|^2) |\Re|.
\end{eqnarray*}
%
%Applying Cauchy-Schwartz inequality, we obtain
% & \le
% C_L N^2 |\int_{t}^T (1+|X_s|^2 + |\bar U_s|^2 + |\bar V_s|^2)^
%Under $\HYP{f}$, we obtain
%which leads to
%(\int_{t}^T |f^{\dot a_s}(X_{s},U^{\dot a}_s, V^{\dot a}_s)
%- f^{a_s}(X_{s}, U^{\dot a}_s, V^{\dot a}_s)|^2 \,\ud s)^{h
% & \le
% C^h_L | N^{\dot a}|^h \sup_{s \in[0,T]} (1+|X_s|^h + | U^{\dot
%a}_s|^h) |\Re|^{h\over2}
%We define
Plugging the last inequality and~(\ref{eqcontrolbA-A1}) in (\ref
{eqteubee2}), we deduce
\[
|\dot\Gamma_t - \Gamma_t |^p \le C^p_L\E_t\Bigl[|N^{\dot
a}|^p\Bigl(\sup_{s \in[0,T]} (1+|X_s|^p + | U_s^{\dot a}|^p)|\Re
|^{p/2}
+ \chi^{|\Re|,p}\Bigr)\Bigr].
\]
Restricting to the case where $t\in\Re$, we deduce from (\ref
{eqcontrolbA-A2}) and~(\ref{eqteubee1}) that
\begin{eqnarray*}
&&
|Y _{t}^{i} - (\tY^{\Re}_{t})^i|^2 \\[-2pt]
&&\qquad\le C^p_L \Bigl( \E_t
\Bigl[|N^{\dot a}|^p\sup_{s \in[0,T]} (1+|X_s|^p + | Y _s|^p)
\Bigr]^{2/p} |\Re|
+ \E_t\bigl[|N^{\dot a}|^p |\chi^{|\Re|,p}\bigr]^{2/p}
\Bigr).
\end{eqnarray*}

Using Cauchy--Schwarz inequality and Proposition~\ref{prstabcont}
with the last inequality, we obtain
\[
|Y _{t}^{i} - (\tY_{t}^{\Re})^i|^2 \le C^p_L \bigl( \beta^X |\Re
| +
\beta^X \mathbb{E}_{t} \bigl[\bigl|\chi^{|\Re|,p}\bigr| ^2
\bigr]^{1/p} \bigr).\vadjust{\goodbreak}
\]

Again using Cauchy--Schwarz inequality and defining $M_{t} := \E
_t
[|\chi^{|\Re|,p}|^2] $, we get
%
%e5.15 ###
%
\begin{equation}\label{4relmiparcourt}
\E\Bigl[\sup_{t\in\Re}|Y ^i_{t} - (\tY^{\Re}_{t})^i|^2\Bigr]
\le
C^p_L \Bigl(|\Re| + \mathbb{E}\Bigl[\sup_{t\in[0,T]} |M_{t}|^{2/ p}
\Bigr]^{1/2}\Bigr) .
\end{equation}

%Using twice Cauchy Schwartz inequality and plugging
%( \E[\sup_{t\in[0,T]} |M^2_{t}|^{2/p}]^{1\over
%2} |\Re|+\E[\sup_{t\in[0,T]} |M^3_{t}|^{2/p}]^{1
%
%where $M^1$, $M^2$ and $M^3$ are the martingales defined by
%M^1_{t} := \E_{t}[ 1 + \sup_{0\le s \le T} |\dot X_s|^{2p} +
%)^{p} ] , \\
%M^2_{t} := \E_t[\sup_{s \in[0,T]} |1+|X_s|^p + | \dot Y_s|^p|^2
%] \mbox{ and }
%M^3_{t} := \E_t[|\chi^{|\Re|,p}|^2] , 0\leq t\leq T
% .

%Combining Burkholder-Davis-Gundy inequality with \eqref{eqcontrolX}
%and Proposition~\ref{4BUShceme}, we get
Combining Burkholder--Davis--Gundy and convexity inequalities with
(\ref{eqcontrolX}), we compute
\begin{eqnarray*}
\E\Bigl[\sup_{t\in[0,T]}|M_{t}|^{2/p}\Bigr]
&\leq& C_{L}^p (|M_{0}|^{2/p}+ \E[|M_{T}|^2
]^{1/p})
\leq C_{L}^p \E\bigl[\bigl|\chi^{|\Re|,p}\bigr|^4\bigr]^{1/p}\\
&\leq& C_{L}^p |\kappa|^{4/p} |\Re|^{2}.
\end{eqnarray*}
Plugging this expression in~(\ref{4relmiparcourt}), we deduce (\ref
{4errordisccontcasgen}) from the condition $\kappa|\Re|\leq L$
and the arbitrariness of $i$.\vspace*{8pt}

\textit{Step} 3. We finally consider the particular case
where the cost functions are constant.
Following the same arguments as in Step 1, we observe that (\ref
{eqcontrolbA-A1}) turns into
\[
\biggl(\int_{t}^T | A^{\dot a}_s -A^a_s |^2 \,\ud s
\biggr)^{p/2} \le
C^p_L |N^{\dot a}|^p \Bigl(1+\sup_{s\in[0,T]}|X_{s}|^p\Bigr)|\Re|^{p/2},
\]
and that $A^{\dot a}_{r}-A^{a}_{r}=0$ for $r\in\Re$.
The same reasoning as in Step 2 then leads to
\[
|Y _{t}^{i} - \tY_{t}^{i}|^2 \le C^2_L \E_t\Bigl[|N^{\dot
a}|^p\sup_{s \in[0,T]} (1+|X_s|^p + | \dot Y_s|^p)\Bigr]^{2/p}
|\Re|.
\]
Using Cauchy--Schwarz and Proposition~\ref{prstabcont} concludes the proof.
%Using Cauchy Schwartz and Burkholder-Davis-Gundy inequalities
%together with \eqref{eqcontrolX} and Proposition~\ref{4BUShceme}
%concludes the proof.
\end{pf}

%We now give a control of the error on the variable $Z$ in the case
%where the maps $c_{ij}$ are constant.
We now present the main result of this section, which allows us to
control the error between the solutions of the continuously and the
discretely obliquely reflected BSDE at any time between $0$ and $T$.
\begin{Theorem} \label{ThmErrorDiscretize2}
Under \textup{(Hz)--(Hc)}, the following holds:
%[Objectif]
%
\begin{eqnarray*}
&&\sup_{t \in[0,T]}\mathbb{E}[|Y _t-\tY^{\Re}_t|^2+|Y _t-Y^{\Re
}_t|^2]+ \mathbb{E}\biggl[\int_0^T|Z _s- Z^{\Re}_s|^2 \,\ud s\biggr]\\
&&\qquad \le% C_
C^\varepsilon_{L}|\Re|^{1/2 - \varepsilon} ,\qquad \varepsilon>0.\nonumber
\end{eqnarray*}
%
%Moreover, if \HYP{f} holds
If, furthermore, the cost functions are constant, the previous estimate
holds true for $\varepsilon= 0$.
\end{Theorem}

In order to prove this theorem, we first state the following lemma
discussing the regularity of $K $.
\begin{Lemma}\label{lereguK}
Under \textup{(Hz)--(Hc)}, there exists some positive process $\eta$
satisfying $\|\eta\|_{\cH^{2}} \le C_{L}$ and such that, for all $i
\in
\cI$,
$\ud K ^{i}_{s} \le\eta_{s} \,\ud s$
in the sense of random measure.
\end{Lemma}
\begin{pf}
We follow here the main idea of the proof of Proposition 4.2 in~\cite
{elkkapparpenque97} and divide the proof in three steps.\vspace*{8pt}

\textit{Step} 1. Fix $i,j \in\cI$. We first observe using It\^{o}--Tanaka
formula, that, under~(Hc),
\[
c^{ij}(X_{t}) = c^{ij}(X_{0}) + \int_{0}^{t}b^{ij}_{s} \,\ud s + \int
_{0}^{t}\nu^{ij}_{s} \,\ud W_{s} - \int_{0}^{t}\ud\Delta^{ij}_{s}
,\qquad 0\le t \le T,
\]
where $\Delta^{ij}$ is an increasing process and
%
%e5.16 ###
%
\begin{equation}\label{eqmajobij}
\|b^{ij}\|_{\cH^{2}} + \|\nu^{ij}\|_{\cH^{2}} \le C_{L} .
\end{equation}

We then introduce $\Gamma^{ij} := Y ^{i} - Y ^{j} +c^{ij}(X) \ge0$.
%% . \label{eqgamma-ij}
Using once again It\^{o}--Tanaka formula, we compute
\begin{eqnarray*}
[\Gamma^{ij}_{t}]^{+} &=& [\Gamma^{ij}_{0}]^{+} + \int
_{0}^{t}\bigl(-f^{i}(X_{s},Y ^{i}_{s},Z ^{i}_{s})+f^{j}(X_{s},Y
^{j}_{s},Z ^{j}_{s}) + b^{ij}_{s}\bigr)\mathbf{1}_{ \{ \Gamma^{ij}_{s}>0
\} } \,\ud s
\\
&&{}+ \int_{0}^{t}(\nu^{ij}_{s}+Z ^{i}_{s}-Z ^{j}_{s})\mathbf{1}_{ \{
\Gamma^{ij}_{s}>0 \} } \,\ud W_{s}\\
&&{}+ \int_{0}^{t}\mathbf{1}_{ \{ \Gamma^{ij}_{s}>0 \} } (- \ud K
^{i}_{s} + \ud K
^{j}_{s} - \ud\Delta^{ij}_{s} )
+\frac12 \int_{0}^{t}\ud L^{ij}_{s}
\end{eqnarray*}
for $0\le t \le T$, where $L^{ij}$ is the local time at $0$ of the
continuous semi-martingale~$\Gamma^{ij}$. Since $\Gamma^{ij}\ge0$ and
$\Delta^{ij}$, $L^{ij}$ are increasing processes, we compute
%
%e5.17 ###
%
\begin{eqnarray} \label{eqtempdKi1}
\mathbf{1}_{ \{ \Gamma^{ij}_{s}=0 \} } \,\ud K ^{i}_{s}
&\le&\bigl(-f^{i}(X_{s},Y ^{i}_{s},Z ^{i}_{s}) + f^{j}(X_{s},Y
^{j}_{s},Z ^{j}_{s}) + b^{ij}_{s} \bigr)\mathbf{1}_{ \{ \Gamma^{ij}_{s}=0
\} } \,\ud s\nonumber\\
&&{}+ \mathbf{1}_{ \{ \Gamma^{ij}_{s}=0 \} } \,\ud K ^{j}_{s}
\nonumber\\[-8pt]\\[-8pt]
&\le& C_{L}\Bigl(1+ |X_{s}|+ \sup_{\ell\in\cI}|Y ^{\ell}_{s}| + \sup
_{\ell
,k \in\cI}|b^{\ell k}_{s}|\Bigr) \,\ud s\nonumber\\
&&{} + \mathbf{1}_{ \{ \Gamma
^{ij}_{s}=0 \} } \,\ud K
^{j}_{s} \nonumber
\end{eqnarray}
for $0\le s \le T$, where we used (Hz) in order to obtain the last
inequality.\vspace*{8pt}

\textit{Step} 2. We now prove that
%
%e5.18 ###
%
\begin{equation}\label{eqdKj=0}
\mathbf{1}_{ \{ \Gamma^{ij}_{s}=0 \} } \,\ud K ^{j}_{s} = 0
\end{equation}
in the sense\vspace*{-2pt} of random measure. We first observe that
$
\mathbf{1}_{ \{ \Gamma^{ij}_{s}=0 \} } \,\ud K ^{j}_{s} = \gamma
^{ij}_{s}\,\ud K
^{j}_{s}$ with $\gamma^{ij}_{s}:=\mathbf{1}_{ \{ \Gamma^{ij}_{s}=0
\} } \mathbf{1}_{ \{ Y^{j}_{s} - \cP^{j}(X_{s},Y _{s}) = 0 \} } $.\vadjust{\goodbreak}
Indeed,\vspace*{2pt} if $\mathbf{1}_{ \{ Y^{j}_{s} - \cP^{j}(X_{s},Y _{s}) > 0 \}
} \,\ud K ^{j}_s$
were a positive random measure on $[0,T]$, this would contradict the
minimality condition~(\ref{eqBSDECOR}) for $K $.

Suppose the existence of a stopping time $\tau$ smaller than $T$, such that
%
%e5.19 ###
%
\begin{equation} \label{eqcontra0}
\Gamma^{ij}_{\tau}=0 \quad\mbox{and}\quad Y ^{j}_{\tau} - \cP
^{j}(X_{\tau},Y
_{\tau}) = 0.
\end{equation}
By definition of the projection $\cP$, we have
%
%e5.20 ###
%
\begin{equation}\label{eqcontra1}
Y ^{j}_{\tau} - \cP^{j}(X_{\tau},Y _{\tau}) = Y ^{j}_{\tau}
- Y
^{k_{\tau}}_{\tau} + c^{j k_{\tau}}(X_{\tau}),
\end{equation}
where\vspace*{1pt} $k_{\tau}$ takes value in $\cI$. Moreover, $Y^{i}_{\tau
}-Y^{k_{\tau}}_{\tau} + c^{ik_{\tau}}(X_{\tau}) \ge0$, which leads,
combined with~(\ref{eqcontra0}) and~(\ref{eqcontra1}), to
$
c^{ij}(X_{\tau})+c^{j k_{\tau}}(X_{\tau}) - c^{ik_{\tau}}(X_{\tau
}) \le0
$
and then contradicts~(\ref{CostStructure}).\vspace*{1pt}

Thus, $\gamma^{ij}_\tau=0$ for any stopping time $\tau$ smaller than
$T$ and we deduce that~$\gamma^{ij}$ is undistinguishable from $0$,
which proves~(\ref{eqdKj=0}).\vspace*{8pt}

\textit{Step} 3. To conclude, using once again the minimality condition
for $K$ in~(\ref{eqBSDECOR}), observe that
$
\ud K ^{i}_{s} = \sum_{j} \mathbf{1}_{ \{ \Gamma^{ij}_{s}=0 \} } \hspace*{-0.2pt}\,\ud
K ^{i}_{s}
\le
\eta_{s} \,\ud s,
$
with $\eta:= C_{L} (1+ |X|+ \sup_{\ell\in\cI}|Y ^{\ell}| + {\sup
_{\ell,k \in\cI}}|b^{\ell k}|) $ which satisfies $\|\eta\|_{\cH
^{2}}\le
C_{L}$, recalling~(\ref{eqcontrolX}),~(\ref{eqestimateia}) and
(\ref{eqmajobij}).
\end{pf}
\begin{pf*}{Proof of Theorem~\ref{ThmErrorDiscretize2}}
Fix $t\in[0,T]$ and introduce $\dotY:= Y - \tY^{\Re}$, $\doY:=
Y -
Y^{\Re}$, $\doZ:= Z -Z^{\Re}$ and $\dof:= f(X,Y ,Z )-
f(X,\tY
^{\Re},Z^{\Re})$.
Applying It\^o's formula to the \textit{c\`{a}dl\`{a}g} process
$|\dotY|^2$, we get
%
%e5.21 ###
%
\begin{eqnarray}
&&|\dotY_t|^2 + \int_t^T|\doZ_s|^2 \,\ud s \nonumber\\[-8pt]\\[-8pt]
&&\qquad= |\dotY_T|^2 - 2\int
_{(t,T]}\dotY_{s-}\,\ud\dotY_s - \sum_{t< s \le T}|\dotY_s -\doY
_{s}|^2.\nonumber
\end{eqnarray}
%
%|\dotY_t|^2 + \int_t^T|\doZ_s|^2 \,\ud s& =& |\dotY_T|^2 - 2
%& & -\sum_{t< s \le T}(|\dotY_s|^2-|\dotY_{s-}|^2-2 \dotY_{s-}.\Delta
%Recall $\dotY_{s-}$ $=$ $\doY_{s}$ for $s\in[0,T]$. Since
%} \int_{(t,T]} \doY_s.\ud K _s \ge0 ,
Recalling that $\dotY_{s-}=\doY_{s}$, $\int_{(t,T]} \doY_s\,\ud K
^{\Re}_s \ge0$ and the Lipschitz property of~$f$, standard arguments
lead to
%2 \int_t^T\dotY_s\dof_s \,\ud s + 2 \int_{t}^T \doY_s.\,\ud K _s} .
%Using standard arguments, we then compute
%
%e5.22 ###
%
\begin{eqnarray}\label{eqcontroldYdZref0}
\mathbb{E}\biggl[|\dotY_t|^2+ \int_t^T|\doZ_s|^2 \,\ud s\biggr] &\le&
C_{L}\mathbb{E}\biggl[\int_t^T \doY_s\,\ud K _s
\biggr]\nonumber\\[-8pt]\\[-8pt]
&\le& C_{L} \sum_{j<\kappa} \mathbb{E}\biggl[\int_{\rj}^{\rjp} \doY
_{s} \,\ud K _{s}\biggr].\nonumber
\end{eqnarray}
%
%recall $\Re$ $=$ $\set{r_{0},\ldots,r_{\kappa}}$.
%with $s+ := \inf\set{t>s| t\in\Re\cup\set{+\infty} } \wedge T$,
%for $s\in[0,T]$.
%

%In the next step, we shall prove that
%in the sense of random measure, for some positive $\eta\in\cH^{2}$.

Using the expression of $\doY$ and Lemma~\ref{lereguK}, we obtain
\[
% \label{eqecrituredY}
\doY_{s} \le\doY_{\rjp} + \int_{s}^{\rjp}(\dof_{u} + \eta
_{u})\,\ud u
-\int_{s}^{\rjp} \doZ_{u} \,\ud W_{u} ,\qquad \rj\le s < \rjp,
j< \kappa.
\]
Combining (Hz),~(\ref{eqcontrolX}),~(\ref{controlContRef}),
(\ref{eqestimateia}) and the fact that $\|\eta\|_{\cH^{2}}\le
C_{L}$, we deduce
\begin{eqnarray*}
\sum_{j<\kappa} \mathbb{E}\biggl[\int_{\rj}^{\rjp} \doY_{s} \,\ud K _{s}\biggr]
& \le &\mathbb{E}\biggl[\sum_{j<\kappa} \int_{\rj}^{\rjp} \!\!\int
_{s}^{\rjp}(\dof_{u} + \eta_{u})\,\ud u \,\ud K _{s}\biggr] \\
&&{}+ \mathbb
{E}\biggl[\sum_{j<\kappa}\int_{\rj}^{\rjp}\doY_{\rjp}\,\ud K _{s}\biggr]
%& \le C_{L}(|\Re| + \esp{\sum_{j<\kappa}\int_{\rj}^{\rjp}\doY_{\rjp}
\\
& \le& C_{L}|\Re| + \mathbb{E}\Bigl[K _{T}\sup_{r \in\Re}|\doY_{r}|\Bigr].
\end{eqnarray*}

Plugging this expression in~(\ref{eqcontroldYdZref0}) and using %
%& \le C_{L}(|\Re| + \esp{K _{T}\sup_{r \in\Re}|\doY_{r}|})
%The proof of the Proposition is then concluded combining the last
%inequality
Cauchy--Schwarz inequality together with~(\ref{4errordisccontcasgen}) and
Proposition~\ref{estprioedsrdrabs} concludes the proof.
\end{pf*}
\subsection{Convergence of the discrete-time scheme}\label{SubSecConvergence}

Combining the previous results with the control of the error between
the discrete-time scheme and the discretely obliquely reflected BSDE
derived in Section~\ref{SubSecScheme}, we obtain the convergence of the
discrete time scheme to the solution of the continuously obliquely
reflected BSDE. In the next theorem, we detail the corresponding
approximation error for different optimal choices of reflection time
step $|\Re|$ with respect to the discrete time step~$|\pi|$.
%Setting $\Re= \pi$, and combining Theorem~\ref{THcvSH} withTheorem
%we obtain the following results.
%
\begin{Theorem}
The following hold:

\begin{longlist}
\item
If \textup{(Hf)--(Hc)} holds, taking $ |\Re| \sim\frac{\log
L_\cP
}{- \eps\log|\pi|}$ for $\eps>0$, we have
\begin{eqnarray*}
&&\sup_{t\in[0,T]} \mathbb{E}[|Y _t - \tYp_t|^2 + |Y _t - \Yp_t|^2] +
\mathbb{E}\biggl[\int^T_{0}|Z _s - \bZp_{s} |^2 \,\ud s\biggr]\\
&&\qquad \le\frac
{C_{L}^\eps
}{[-{\log}
(|\pi|)]^{1/2-\eps}}.
\end{eqnarray*}

\item If $f$ does not depend on $z$ and $|\pi|L<1$, taking
similar grids $\Re=\pi$, we have
\[
\sup_{i\le n}\mathbb{E}[|Y _\ti- \Yp_\ti|^2 + |Y _\ti- \tYp_\ti
|^2 ]
\le C_L^\eps|\pi|^{1-\eps} ,\qquad \eps>0.
\]
Moreover, under \textup{(Hc)},
\[
\sup_{t \in[0,T]} \mathbb{E}[|Y _t - \Yp_t|^2 + |Y _t - \tYp_t|^2]
\le
C^{\varepsilon}_L |\pi|^{1/2 - \eps} ,\qquad \eps>0.
\]

\item Under \textup{(Hc)}, if $f$ does not depend on $z$ and $|\pi|L<
1$, taking $|\Re|\sim|\pi|^{2/3}$, we get
\[
\mathbb{E}\biggl[\int^T_{0}|Z _s - \bZp_{s} |^2 \,\ud s\biggr] \le C_{L}^{\eps}
|\pi
|^{1/3-\eps},\qquad \eps>0.
\]

\item Furthermore, for constant cost functions, \textup{(ii)} and
\textup{(iii)} hold true with $\eps= 0$.\vadjust{\goodbreak}
\end{longlist}
\end{Theorem}
%
%Under \HYP{f}, the following holds %discrete-time scheme $(\tYp,\bZp)$
%converges to $(\tY,Z)$:
%for all $\eps>0$ and $|\pi|$ small enough.
%
\begin{pf}
For $\eps>0$, setting $\Re$ such that $ |\Re| \sim\frac
{\log
L_\cP}{- \eps\log|\pi|}$, we obtain, combining Proposition \ref
{ThmErrorScheme} and Theorem~\ref{ThmErrorDiscretize2}, that
\begin{eqnarray*}
&&\sup_{t\in[0,T]} \mathbb{E}[|Y _t - \tYp_t|^2 + |Y _t - \Yp_t|^2]
+ \mathbb{E}\biggl[\int^T|Z _s - \bZp_{s} |^2 \,\ud s\biggr]\\
&&\qquad\le
C^{\eps}_{L} \biggl[ \biggl(\frac{-1}{\log(|\pi|)}\biggr)^{1/2
-\eps} \vee|\pi|^{1/2 - \eps} \biggr].
\end{eqnarray*}
Therefore, (i) is proved. Furthermore, (ii), (iii) and (iv) are direct
consequences of Theorems~\ref{THcvSH} and~\ref{ThmErrorDiscretize} or
\ref{ThmErrorDiscretize2}.
%The proof is concluded setting $\eta= \frac{\eps}{\log L}$.%$\eta=
\end{pf}

%laisser avec des $\varepsilon$ dans le theoreme ou optimiser la borne...,
%c'est le log qui ralentit qd $\pi$ est tres fine...ctb}

%This convergence result is improved when $f$ does not depend on $z$.
%Setting $\Re= \pi$, and combining Theorem~\ref{THcvSH} with Theorem
%we obtain the following results.

% (i) If $f$ does not depend on $z$ and $|\pi|L$ $<$ $1$, the following
%holds
%and
%& C^{\varepsilon}_L |\pi|^{ \frac12 - \eps} ,
%for all $\varepsilon> 0$.

%
%(ii) If the cost functions are constant, the previous inequality holds
%true with $\varepsilon= 0$.

%We now turn to the approximation of the $z$ component. Setting $\Re$,
%such that $|\Re| \sim|\pi|^{\frac13}$, we obtain combining
%(i) If $f$ does not depend on $z$ and $|\pi|L$ $<$ $1$, the following
%holds
%for $\varepsilon>0$.

%(ii) If the cost functions are constant, the previous inequality holds
%true with $\varepsilon= 0$.

%%%%%%%%%%%%%%%%%%%%%%
%%%%%%%%%%%%%%%%%%%%%%
%%%%%%%%%%%%%%%%%%%%%%
%%%%%%%%%%%%%%%%%%%%%%

%!TEX root = main.tex

%
\begin{appendix}\label{app}
\section*{Appendix}
%s5.4 ###
\subsection{A priori estimates for discretely RBSDEs}
We collect here the proofs for a priori estimates given in Propositions
\ref{estprioedsrdrabs} and~\ref{corhatA}.
\begin{pf*}{Proof of Proposition~\ref{estprioedsrdrabs}}
Observing that on each interval $[\rj,\rjp)$, $(Y^{\Re},\tY^{\Re
},Z^{\Re
})$ solves
a standard BSDE, existence and uniqueness follow from a
concatenation procedure and~\cite{parpen90}.
The rest of the proof divides in two steps controlling separately $\tY
^{\Re}$ and $(Z^{\Re},K^{\Re})$.\vspace*{8pt}

\textit{Step} 1. \textit{Control of $\tY^{\Re}$.}

As in the proof of Theorem 2.4 in~\cite{hz08}, we consider two
nonreflected BSDEs
bounding $\tY^{\Re}$.

Define the $\R^d$-valued random variable $\breve\xi$ and the
random map $\breve F$ by $(\breve\xi)^j:=\sum_{i=1}^d |\xi|^i$ and
$(\breve F)^j := \sum_{i=1}^d |(F)^i|$ for $1\le j\le d$.\vspace*{1pt}

We then denote by $(\breve Y, \breve Z)\in(\Sc^2\times\cH^2)^\Ic$ the
solution to the following nonreflected BSDE:
%
%e5.23 ###
%
\setcounter{equation}{0}
\begin{equation}\label{4BSDEbreve}
\breve Y_{t} = \breve\xi+ \int_{t}^T\breve F(s,\breve
Y_{s},\breve Z_{s})\,\ud s-\int_{t}^T\breve Z_{s} \,\ud W_{s} ,\qquad
0\leq t \leq T.
\end{equation}
Since all the components of $\breve Y$ are similar, $\breve Y\in\cC$.

We also introduce $({}^0Y,{}^0Z)$ the solution to
the BSDE
\[
\label{4BSDEZERO}%\label{edsrsc}
{}^0Y_{t} = \xi+\int_{t}^T F({s},{}^0Y_{s},{}^0 Z_{s})\,\ud s-\int
_{t}^T{}^0 Z_{s} \,\ud W_{s} ,\qquad 0\leq t\leq T.
\]
Using\vspace*{1pt} a comparison argument on each interval $[\rj,\rjp)$ and the
monotony property of $\cP$, we straightforwardly deduce
${}^0Y \preceq Y^{\Re} \preceq\breve Y$.\vspace*{1pt}

Since $({}^0Y, \breve Y)$ are solutions to standard nonreflected BSDEs,
usual arguments lead to
%
%e5.24 ###
%
\begin{equation} \label{eqespsupyp}
\sup_{0\le s\le T}|\tY^{\Re}_s|^p \le\sup_{0\le s\le T}|{}^0Y_s|^p
+ \sup_{0\le s\le T}|\breve Y_s|^p =: \bar{\beta},
\end{equation}
where the positive random variable $\bar{\beta}$ satisfies classically
$\mathbb{E}[\bar{\beta}] \le C_{L}$, under condition
(\ref{Cp}) for a
given $p\ge2$.\vspace*{8pt}

\textit{Step} 2. \textit{Control of $(Z^{\Re},K^{\Re})$.}

We fix $t\le T$ and applying It\^o's formula to the \textit{c\`{a}dl\`
{a}g} process
$|\tY^{\Re}|^2$ on $[0,t]$ to derive
\begin{eqnarray*}
|\tY^{\Re}_t|^2 &=& |\tY^{\Re}_0|^2 + 2 \int_{(0,t]}\tY^{\Re
}_{s-}\,\ud
\tY^{\Re}_s + \int_{(0,t]}|Z^{\Re}_s|^2 \,\ud s\\
&&{}+ \sum_{s\le t}(|\tY^{\Re}_s|^2-|\tY^{\Re}_{s-}|^2- 2 \tY^{\Re}_{s-}
\Delta Y ^{\Re}_s).
\end{eqnarray*}
Since the last term on the right-hand side is nonnegative, we deduce
that % from \reff{eq def Yg}
\begin{eqnarray*}
|\tY^{\Re}_t|^2 + \int_t^T|Z^{\Re}_s|^2 \,\ud s &\le&
|\tY^{\Re}_T|^2 + 2 \int_t^T \tY^{\Re
}_{s-}F(s,\tY
^{\Re}_s,Z^{\Re}_s)\,\ud s\\
&&{}+ 2 \int_{(t,T]}\tY^{\Re}_{s-} \,\ud K ^{\Re}_s + 2 \int_t^T (
Z^{\Re
}_s\tY^{\Re}_{s}) \,\ud W_s.
\end{eqnarray*}
Using standard arguments, together with
(\ref{eqespsupyp}) and~(\ref{Cp}) for a fixed $p\ge2$, we compute
%
%e5.25 ###
%
\begin{equation} \label{eqapco2c}\quad
\int_t^T|Z^{\Re}_s|^2\,\ud s \le C_{L} \biggl( \bar{\beta}^{2
/
p} + \bar{\beta}^{1/p}(K ^{\Re}_T-K ^{\Re}_t) +
\int_t^T ( Z^{\Re}_s\tY^{\Re}_s) \,\ud W_s \biggr).
\end{equation}

Moreover, we get from~(\ref{eqBSDEDORabs}) and~(\ref{Cp}) that
%
%e5.26 ###
%
\begin{equation}\label{eqapco2cbis}
|K ^{\Re}_T - K ^{\Re}_t|^{2} \le
C_{L} \biggl[\bar{\beta}^{2/p} + \int_t^T|Z^{\Re}_s|^2 \,\ud s +
\biggl(\int_t^T Z^{\Re}_s \,\ud W_s\biggr)^2 \biggr].
\end{equation}

Combining~(\ref{eqapco2c}) and~(\ref{eqapco2cbis}) we obtain
%
%e5.27 ###
%
\begin{eqnarray}\label{eqapco2d}
\int_t^T|Z^{\Re}_s|^2\,\ud s &\le&\frac{C_{L}}{\eps} \bar{\beta
}^{2/p} + \eps
\int_t^T|Z^{\Re}_s|^2 \,\ud s + \eps\biggl(\int_t^TZ^{\Re}_s \,\ud W_s
\biggr)^2 \nonumber\\[-8pt]\\[-8pt]
&&{}+C_{L} \int_t^T ( Z^{\Re}_s\tY^{\Re}_s)\,\ud W_s\nonumber
\end{eqnarray}
for any $\eps>0$.
Elevating the previous estimate to the power $p/2$, it follows from
Burkholder--Davis--Gundy inequality that
\begin{eqnarray*}
&&\mathbb{E}_{t} \biggl[\biggl(\int_t^T|Z^{\Re}_s|^2\,\ud s
\biggr)^{p/2}\biggr]\\
&&\qquad\le
C^{p}_L \biggl( \eps^{-{p/2}} \mathbb{E}_{t} [\bar
{\beta}] + \eps
^{p/2}\mathbb{E}_{t} \biggl[\biggl(\int_t^T|Z^{\Re}_s|^2\,\ud s\biggr)^
{p/2}\biggr] \\
&&\hspace*{109.5pt}{} + \mathbb{E}_{t} \biggl[\biggl(\int_t^T |Z^{\Re}_s\tY
^{\Re}_s |^2 \,\ud s\biggr)^{p/4}\biggr]\biggr)\\
&&\qquad \le  C^{p}_L \biggl(\eps^{-{p/2}}\mathbb{E}_{t} [\bar
{\beta}] + \eps
^{-{p/2}}\mathbb{E}_{t} \Bigl[\sup_{s\in[t,T]}|\tY^{\Re
}_s|^{p}\Bigr] \\
&&\hspace*{81.7pt}{} + \eps
^{p/2}\mathbb{E}_{t} \biggl[\biggl(\int_t^T|Z^{\Re}_s|^2\,\ud s\biggr)^
{p/2}\biggr] \biggr).
\end{eqnarray*}
Using~(\ref{eqespsupyp}) and~(\ref{Cp}), we deduce, for
$\eps$
small enough,
%
%e5.28 ###
%
\begin{equation} \label{eqcontrolintZt-T}
\mathbb{E}_{t} \biggl[\biggl(\int_t^T|Z^{\Re}_s|^2\,\ud s\biggr)^
{p/2}\biggr] \le C_L^p\mathbb{E}_{t} [\bar{\beta}
] .
\end{equation}
Taking\vspace*{1pt}~(\ref{eqapco2cbis}) up to the power $\frac{{p}}{2}$, and
combining Burkholder--Davis--Gundy inequality with~(\ref{eqcontrolintZt-T})
yields $\mathbb{E}_{t} [|K ^{\Re}_T-K ^{\Re}_t|^p]
\le C_L^p\mathbb{E}_{t} [\bar{\beta}]$, which
concludes\vspace*{1pt} the proof of the proposition, recalling~(\ref{Cp}).
\end{pf*}
\begin{pf*}{Proof of Proposition~\ref{corhatA}}
Fix $(t,i)\in[0,T]\times\Ic$ and $p\ge2$. According to the
identification of $(U^{a^*},V^{a^*})$ with
$(\tY^{a^*},Z^{a^*})$, obtained in the proof of Theorem~\ref{THswR},
we deduce from Proposition~\ref{estprioedsrdrabs} the expected
controls on $U^{a^*}$ and $V^{a^*}$.
Writing the equation satisfied by $(U^{ a^*},V^{ a^*})$ and using
standard arguments for BSDEs, we observe that
\[
\mathbb{E}_{t} [|A_{T}^{a^*}|^p] \le C^p_{L} \biggl(
\mathbb{E}_{t} \biggl[\sup_{s\in[t,T]}|U^{a^*}_{s}|^p +
\biggl(\int_{t}^T|V^{a^*}_{s}|^2\,\ud s \biggr)^{p/2}\biggr] +
|A^{a^*}_{t}|^p \biggr).
\]
By definition of $a^*$ and~(\ref{CostStructureabs}), we have
$|A_{t}^{a^*}|\leq\max_{k\neq i}|C^{i,k}_{t}|$, which plugged in
the previous inequality leads to $ \mathbb{E}_{t}
[|A_{T}^{a^*}|^p] \le
C^p_{L} \mathbb{E}_{t} [\bar{\beta}]$, recalling
(\ref{Cp}).

We finally complete the proof, noticing from~(\ref{CostStructureabs})
that $\mathbb{E}_{t} [|N^{a^*}|^p] \leq C_{L}^p\times
\mathbb{E}_{t} [|A_{T}^{a^*}|^p]$.
% which gives the waited result.
% we obtain as in the proof of Theorem~\ref{THswR} (ii) that
% \esp{|A_{T}^{a^*}|^p} \le C_L^p .
\end{pf*}

%s5.5 ###
\subsection{A priori estimates for the Euler scheme}\label{secap2}

This paragraph provides the proof of Proposition~\ref{4BUShceme},
concerning a priori estimates for the Euler scheme associated to RBSDEs.
\begin{pf*}{Proof of Proposition~\ref{4BUShceme}}
The proof follows exactly the same arguments as in Step 1 of the proof
of Proposition~\ref{estprioedsrdrabs} above. The only difficulty is
the use of a comparison argument for Euler scheme that we provide below
in Lemma~\ref{PropCompaScheme}.
\end{pf*}

We detail here a comparison theorem for discrete-time schemes of BSDEs
in the case where the driver does not depend on the variable $z$.

For $k=1,2$, let $\xi_{k}$ be a square integrable random variable and
$\psi_{k} \dvtx\R^m\times\R^d\rightarrow\R$ an $L$-Lipschitz generator
function.\vadjust{\goodbreak} We suppose that $\xi_{1}\geq\xi_{2}$ and $\psi_{1}
\geq\psi_{2}$ on $\R^m\times\R^d$. For a time grid $\pi$, we denote
by $Y^{\pi,k}$ %$(Y^{\pi,k},\bar Z^{\pi,k})$, $k=1,2$
the discrete-time scheme starting from the terminal condition
$
Y^{\pi,k}_T := \xi_{k}
$
and computing recursively, for $i = n-1,\ldots,0$,
%
%e5.29 ###
%
\begin{equation}\label{schemethmcomp}
% \{
% \begin{array}{rcl}
% \bar Z ^{\pi,k}_{\ti} &=& (\tip-\ti)^{-1}\EF{ Y^{\pi,k}_\tip(W_\tip-W_
% \\
Y^{\pi,k}_{\ti} = \mathbb{E}[Y^{\pi,k}_\tip\mid\mathcal
{F}_{\ti}] + (\tip- \ti)
\psi
_{k}(\Xp_\ti,Y^{\pi,k}_\ti).
% \\
% \Yp_{\ti}&=& \tYp_\ti\Ind{\ti\notin\Re} + \cP(\Xp_\ti,\tYp_\ti)
% \end{array}
%
\end{equation}

\setcounter{Lemma}{0}
\begin{Lemma}\label{PropCompaScheme}
%Let $\xi_{1}$ and $\xi_{2}$ be two square integrable random variables
%and $f_{1}$ and $f_{2}$ two generators. We suppose that $\xi_{1}$ $
%denote by $(Y^{\pi,i},\bar Z^{\pi,i})$, $i=1,2$, the discrete-time
%scheme
% METTRE LE DISCRETE TIME SCHEME
For any $\pi$ such that $|\pi| L<1$, we have $Y^{\pi,1}_{\ti}
\geq
Y^{\pi,2}_{\ti}$, $i\le n$.
\end{Lemma}
\begin{pf}
%Recall that$\pi= \set{t_{0}=0,\ldots,t_{n}=T}$.
Since the result holds true on the grid point $t_n=T$ and follows from
a backward induction on $\pi$, we just prove $Y^{\pi,1}_{t_{n-1}}
\geq Y^{\pi,2}_{t_{n-1}}$.
%Notice that it suffices to prove $Y^{\pi,1}_{t}$ $\geq$ $Y^{
%in a backward manner at each step of $\pi$.
Using~(\ref{schemethmcomp}), we compute
%For the definition of the scheme we have
%
%e5.30 ###
%
\begin{eqnarray}\label{4deltaYappendix}\qquad
Y^{\pi,1}_{t_{n-1}}-Y^{\pi,2}_{t_{n-1}} &=& \E_{t_{n-1}} [\xi
_{1}-\xi
_{2}\mid\Fc_{t_{n-1}}] + (T-{t_{n-1}})\Lambda_{n-1}
(Y^{\pi
,1}_{t_{n-1}}-Y^{\pi,2}_{t_{n-1}})\nonumber\\[-8pt]\\[-8pt]
&&{} + \Delta_{n-1},\nonumber
\end{eqnarray}
where $\Delta_{n-1} := \psi_{1}(X^\pi_{t_{n-1}}Y^{\pi
,2}_{t_{n-1}})-\psi
_{2}(X^\pi_{t_{n-1}}Y^{\pi,2}_{t_{n-1}})\geq0$ and
%
%e5.31 ###
%
\begin{equation}\qquad
\Lambda_{n-1} :=
\cases{
\displaystyle \frac{\psi_{1}(X^\pi_{t_{n-1}}Y^{\pi,1}_{t_{n-1}})-\psi_{1}(X^\pi
_{t_{n-1}}Y^{\pi,2}_{t_{n-1}})}{Y^{\pi,1}_{t_{n-1}}-Y^{\pi
,2}_{t_{n-1}}}, &\quad if $Y^{\pi,1}_{t_{n-1}}-Y^{\pi,2}_{t_{n-1}} \neq0$,
\vspace*{2pt}\cr
0, &\quad else.}
\end{equation}
%
%and
Since $\psi_{1}$ is $L$-Lipschitz, the condition $|\pi|L<1$,
implies $(T-t_{n-1})\Lambda_{n-1}<1$. Plugging this estimate,
$\Delta_{n-1} \geq0$ and $\xi_{1}\geq\xi_{2}$ and $\psi_{1}$ in
(\ref{4deltaYappendix}), the proof is complete.
\end{pf}

%
%%%%%%%%%%%%%%%%%%%%%%%%%%%%%%
%

%s5.6 ###
\subsection{A priori estimates for continuously RBSDEs}\label{secap3}
This last paragraph is dedicated to the proof of Proposition~\ref{prstabcont}.
\begin{pf*}{Proof of Proposition~\ref{prstabcont}}
The proof of~(\ref{eqestimateia}) is a direct adaptation of the
proof of Proposition~\ref{estprioedsrdrabs}.
The only difference is in Step 1: we approximate $(Y ,Z ,K )$ by a
sequence of penalized BSDEs (see the proof of Theorem 2.4 in \cite
{hz08} or Step 3 in the proof of Theorem~\ref{ThmErrorDiscretize2}) which
are bounded by ${}^0Y$ and~$\breve Y$. Estimate~(\ref{eqestimateib})
follows from the exact same arguments as the one used in the proof of
Proposition~\ref{corhatA}.
%
%(ii) To prove the stability property, we follow line by line the Step
%1 and Step 2 in proof of Proposition~\ref{prstabY}. The only
%difference is that the supremum is taken over $[0,T]$ instead of $
\end{pf*}
\end{appendix}

%suskaldyti doi

% imsref loaded by lrinkeviciute, 2011-08-08 10:29:32
%

%
\printaddresses

\end{document}